\documentclass[11pt]{article}

\usepackage{graphicx}
\usepackage{wrapfig}
\usepackage{substr}
\usepackage[tight]{minitoc}
\usepackage{subfigure}
\usepackage{longtable}
\usepackage[all]{xy}
\usepackage{xcite}
\usepackage{xr}

\usepackage{amsfonts}
\usepackage{nicefrac}

\usepackage{amsmath,amssymb,amsthm}
\usepackage{mathtools}

\newcommand{\sR}{{\R^{d}}}

\newcommand{\und}[2]{\underset{#2}{#1}}

\newtheorem{teo}{Theorem}%[section] 
\newtheorem{lem}{Lemma}%[section] 
\newtheorem{assum}{Assumption}%[section] 

\usepackage{fullpage}
\usepackage{cancel}
\usepackage[utf8]{inputenc} % allow utf-8 input
\usepackage[T1]{fontenc}    % use 8-bit T1 fonts
\usepackage{hyperref}       % hyperlinks
\usepackage{url}            % simple URL typesetting
\usepackage{booktabs}       % professional-quality tables
\usepackage{amsfonts}       % blackboard math symbols
\usepackage{nicefrac}       % compact symbols for 1/2, etc.
\usepackage{microtype}      % microtypography

\usepackage{natbib}
\usepackage{sidecap}

\usepackage{amssymb, amsmath, amsthm, latexsym}
\usepackage{url}
\usepackage{algorithm}
\usepackage{algorithmic}
\usepackage{tabularx}
\usepackage{paralist}
\usepackage{mathtools}

\usepackage{bbm} %for indicators 1

\usepackage{makecell}
\usepackage{multirow}
\usepackage{booktabs}

\usepackage{nicefrac}       % nice fractions

\usepackage[flushleft]{threeparttable} % http://ctan.org/pkg/threeparttable

\usepackage{caption}
\usepackage{subcaption}
\usepackage{multirow}
\usepackage{colortbl}
\definecolor{bgcolor}{rgb}{0.8,1,1}
\definecolor{bgcolor2}{rgb}{0.8,1,0.8}
\definecolor{niceblue}{rgb}{0.0,0.19,0.56}

\hypersetup{colorlinks,linkcolor={blue},citecolor={blue},urlcolor={blue}}

\usepackage[dvipsnames]{xcolor}
\usepackage{pifont}

\newcommand{\R}{\mathbb{R}}

\def\<#1,#2>{\left\langle #1,#2\right\rangle}

\usepackage{mdframed} 
\usepackage{thmtools}

\definecolor{shadecolor}{gray}{0.9}
\declaretheoremstyle[
headfont=\normalfont\bfseries,
notefont=\mdseries, notebraces={(}{)},
bodyfont=\normalfont,
postheadspace=0.5em,
spaceabove=1pt,
mdframed={
 skipabove=8pt,
 skipbelow=8pt,
 hidealllines=true,
 backgroundcolor={shadecolor},
 innerleftmargin=4pt,
 innerrightmargin=4pt}
]{shaded}

\declaretheorem[style=shaded,within=section]{definition}

\declaretheorem[style=shaded,sibling=definition]{remark}

\usepackage{xspace}

\usepackage[colorinlistoftodos,bordercolor=orange,backgroundcolor=orange!20,linecolor=orange,textsize=scriptsize]{todonotes}

\usepackage{hyperref}
\graphicspath{{plots/}}

\usepackage{makecell}

\usepackage{accents}
\newlength{\dhatheight}

\usepackage{pgfplotstable} % loads pgfplots which loads tikz
\usetikzlibrary{automata, positioning, arrows, shapes, fit, calc, intersections}
\usepgfplotslibrary{statistics}
\include{figure}

\usepackage{pgfplotstable} % loads pgfplots which loads tikz
\usetikzlibrary{automata, positioning, arrows, shapes, fit, calc, intersections}
\usepgfplotslibrary{statistics}
\include{figure}

\renewcommand{\leq}{\leqslant}

\newcommand\swapifbranches[3]{#1{#3}{#2}}
\makeatletter
\MHInternalSyntaxOn
\patchcmd{\DeclarePairedDelimiter}{\@ifstar}{\swapifbranches\@ifstar}{}{}
\MHInternalSyntaxOff
\makeatother

\DeclarePairedDelimiterX{\inp}[2]{\langle}{\rangle}{#1, #2}
\DeclarePairedDelimiterX{\abs}[1]{\lvert}{\rvert}{#1}
\DeclarePairedDelimiterX{\roundup}[1]{\lceil}{\rceil}{#1}
\DeclarePairedDelimiterX{\cbr}[1]{\{}{\}}{#1} % curly bracket
\DeclarePairedDelimiterX{\rbr}[1]{(}{)}{#1} % round bracket
\DeclarePairedDelimiterX{\sbr}[1]{[}{]}{#1} % 

\title{\bf Computer-Aided Analysis of Gradient Descent under Generalized Smoothness}
\author{Valery Krivchenko \\
Moscow Institute of Physics and Technology\thanks{Moscow Institute of Physics and Technology,
Institutskiy Pereulok, 9, Dolgoprudny, Moscow Oblast, 141701, Russia},\\
\texttt{krivchenko.vo@phystech.edu}}

\date{\today}

\begin{document}

\maketitle

\abstract{Under classical assumptions, performance of first-order methods has been studied quite comprehensively with the developed theoretical apparatus of modern optimization. The development of the Performance Estimation Problems (PEP) framework led to even deeper understanding and allowed computer-assisted derivation of exact worst-case guarantees. Unfortunately, PEP is not directly applicable for systems that involve SDP-incompatible constraints. This includes a wide range of non-standard assumptions, some of them are highly relevant for practice such as generalized smoothness. In this paper proposes an approach of utilizing template solutions to derive guarantees in settings that do not admit an equivalent SDP formulation of exact worst-case performance problem. The approach is showcased by studying the performance of the gradient descent method for convex functions that also satisfy $(L_0,L_1)$-smoothness assumption. Despite the inherent conservatism of the approach, this paper improves and generalizes an existing convergence guarantee to showcase the approach.}


\noindent\textbf{Note:} This is extended version of the paper \cite{Krivchenko2026} published in \textit{Doklady Mathematics} in 2026, doi: 10.7868/S3034504926030011. It is much more detailed overall and contains descriptions of numerical PEP experiments. %\footnote{This is extended version of the paper \cite{Krivchenko2026} published in \textit{Doklady Mathematics} in 2026, doi: 10.7868/S3034504926030011. It is much more detailed overall and contains descriptions of numerical PEP experiments.}

\section{Introduction}

Consider the unconstrained minimization problem
\begin{align*}
    \min_{x \in \R^d} \; f(x),
\end{align*}

where the objective function $f$ satisfies a set of standard assumptions such as convexity, smoothness, and related regularity conditions. In high-dimensional regimes arising, for instance, in machine learning and data analysis, first-order methods constitute a central algorithmic paradigm due to their scalability, low per-iteration cost, and favorable empirical performance \cite{NIPS2007_0d3180d6, Nesterov2013-mm}. Under classical assumptions, most notably convexity and Lipschitz continuity of the gradient, the convergence behavior of these methods is by now well understood \cite{Nesterov2013-mm, beck2017first}.

Over the past decade, this understanding has been significantly deepened through the development of computer-aided analysis techniques, most prominently the Performance Estimation Problem (PEP) framework \cite{drori2014performance}. The PEP framework enables the derivation of worst-case performance guarantees for first-order methods by reformulating their analysis as a convex optimization problem, typically a semidefinite program (SDP). This approach has led to sharp, and in many cases provably optimal, convergence bounds for a broad class of algorithms and problem settings \cite{taylor2017smooth, De_Klerk2017-fm, doi:10.1137/16M108104X, Taylor2018-gk}.

More recently, there has been growing interest in relaxing classical smoothness assumptions in order to better capture the structure of modern optimization problems, in particular those arising in deep learning \cite{ICLR2025_28a9255e}. One prominent example is generalized smoothness, including the $(L_0,L_1)$-smoothness condition, which allows the Lipschitz constant of the gradient to depend on the function value or the norm of the gradient \cite{10.5555/3666122.3667871}. A number of first-order methods have been analyzed under $(L_0,L_1)$-smoothness and related assumptions \cite{10.5555/3495724.3497025, ICLR2025_28a9255e, gorbunov2025methods, lobanov2025linearconvergencerateconvex}. Among the latest works on $(L_0,L_1)$-smoothness, in \cite{lobanov2025linearconvergencerateconvex} there were derived improved complexity bounds for $(L_0,L_1)$-\textit{gradient descent} (GD), \textit{normalized GD} (NGD), \textit{GD with clipping} (Clip-GD) assuming $(L_0,L_1)$-smoothness and convexity. A special attention was paid to special gradient regimes. Thus, $(L_0,L_1)$-GD is guaranteed to converge like $f(x^N)- f(x_*) \leqslant \left(1-\frac{1}{4L_1R}\right)^N F_0$
when $\|\nabla f(x^k) \| \geqslant L_0/L_1$. And for NGD, the guarantee is $f(x^N)- f(x_*) \leqslant \left(1-\frac{\eta}{2R}\right)^N F_0$ for $\eta \leqslant \frac{c}{L_0+cL_1}$ if $\|\nabla f(x^k) \| \geqslant c$ for NGD. 

However, these analyses typically rely on problem-specific and often highly customized analytical arguments. As a result, it remains largely unclear whether the obtained convergence bounds are tight or whether the analyses themselves are optimal. This issue becomes even more pronounced when additional structural assumptions are introduced, such as the presence of different gradient regimes or local smoothness properties.

Unfortunately, for most non-standard function classes—including the $(L_0,L_1)$-smooth class—the PEP framework cannot be applied directly. The fundamental difficulty lies in the lack of necessary and sufficient interpolation conditions that admit a SDP-compatible representation \cite{goujaud2023fundamental}. A naive workaround would be to replace the problematic constraints by their SDP-compatible approximations and solve the resulting relaxed PEP. However, constructing such approximations in a principled way, while retaining meaningful guarantees, is itself a highly nontrivial task.

The goal of this paper is to introduce an approach that enables the effective use of the PEP framework for obtaining worst-case guarantees in settings governed by non-standard assumptions that cannot be characterized by an equivalent system of SDP-compatible constraints. The core idea is to fix constraints by treating certain components of the constraints as fixed parameters. More specifically, the proposed approach consists of three steps. The first step is to derive a set of inequalities that characterize the function class of interest. The second step is to obtain a relaxation of this system by designating certain algebraic components of the inequalities as constant parameters. Because the resulting formulation is a relaxation of the original constraints, its corresponding dual PEP yields valid guarantees. When available, additional information about the dynamics of the algorithm may be incorporated to reduce conservatism. The final step is formulate and solve the resulting PEP, which now takes the form of tractable SDP, and extract the desired performance bounds.

From a different perspective, the final step amounts to instantiating an existing parameterized PEP solution with problem-specific parameters \cite{taylor2017smooth}. While the first two steps necessarily require manual, problem-dependent derivations, the use of PEP can substantially reduce analytical effort once a suitable parameterized SDP formulation is available. Even when such a formulation must be derived from scratch, the PEP framework helps avoid unnecessary conservatism introduced by suboptimal or incomplete analytical arguments.

The main limitation of the proposed approach is that it does not, in general, yield exact worst-case guarantees. Nevertheless, it offers several important advantages. Despite the inherent relaxation, the resulting bounds may still be significantly less conservative than existing alternatives. Moreover, the approach enables a partially automated analysis pipeline, thereby making a concrete step toward closing the gap between modern generalized smoothness assumptions and computer-aided worst-case analysis.

\paragraph{Contributions.} Our contribution can be summarized as follows:

\begin{itemize}\itemsep=-2pt

\item The paper proposes several parameterized SDP programs to be used as templates for studying convergence of gradient descent in generalized smoothness setting. Unlike usual PEP schemes where only the distance to optimum is bounded, both of our template schemes allow arbitrarily set initial gradient norms and functional errors. We derive necessary conditions for program's feasibility and numerically analyze their sufficiency. \\

\item We derive the upper bounds for template scheme's solutions and verify numerically that the acquired bounds are exact. Because different parameter regimes may lead to different closed form solutions, we analyze all the regimes that are relevant to generalized smoothness applications.\\

\item  We obtain equivalent inequalities related to $(L_0,L_1)$-smooth convex functions that do not bound distances $\|x-y\|$. With the inequalities, we are able to apply the proposed PEP template solutions to acquire convergence guarantees for $\|\nabla f(x^k) \| \geqslant L_0/L_1$ regime. Thus we managed to achieve improvement over state-of-the-art convergence guarantees for $(L_0,L_1)$-GD in the setting of $(L_0,L_1)$-smooth convex minimization.

\end{itemize}

\section{Preliminaries}

We work in a finite-dimensional real vector space $\R^d$ equipped with the standard inner product $\langle \cdot, \cdot \rangle$ and the induced Euclidean norm $\|\cdot\|_2$. Throughout the paper, $f: \sR\rightarrow \R$ denotes a proper, closed, and convex function, differentiable on $\sR$.

In this paper we will consider only unconstrained minimization problem 
\begin{align*}
    \min_{x \in \R^d} \; f(x),
\end{align*}
where $f$ may satisfy additional assumptions of standard or generalized smoothness. More formally:

%Assumption: convexity
\begin{assum} [Convexity.] \label{assum:convex} Function $f: \sR\rightarrow \R$ is $\mu$-strongly-convex with $\mu \geqslant 0$:
\begin{align}
    f(y) \geqslant f(x) + \langle f(x), y-x \rangle + \frac{\mu}{2}\|y-x\|_2^2, \hspace{4pt} x,y \in \R^d.
\end{align}
\end{assum}

Next, the standard definition of Lipschitz-smooth gradient:

%Lipshitz smoothness
\begin{assum}[Lipshitz smoothness] \label{assum:smooth} Function $f: \sR\rightarrow \R$ is $L$-smooth, i.e. for $x,y \in \sR$ we have
\begin{align}
    \| \nabla f(x) - \nabla f(y)\|_2 \leqslant L \|x-y\|_2.
\end{align}
\end{assum}

The main focus of the paper is our semi-automated approach of deriving upper bounds for setups that do not have SDP-compatible descriptions. We showcase our method in a $(L_0,L_1)$-generalized smoothness setting that drew a lot of interest due to its flexibility better capacity for modeling loss functions in DL. In literature, there can be found several definitions of $(L_0,L_1)$-smoothness with varying degrees of broadness. In this paper we shall consider $\mathcal L^*_{sym}(1)$ \cite{chen_generalized_smooth_2023} functions that are commonly known as symmetric $(L_0,L_1)$-smooth functions.

%Symmetric smoothness
\begin{assum} [Symmetric $(L_0,L_1)$-smoothness] \label{assum:sym_smooth} Function $f: \sR\rightarrow \R$ is symmetrically $(L_0,L_1)$-smooth, i.e. for $x,y \in \sR$ we have
\begin{align}
    \| \nabla f(x) - \nabla f(y)\|_2 \leqslant \left (L_0+L_1\und{\sup}{u\in[x,y]}\|\nabla f(u)\|_2 \right) \|x-y\|_2.
\end{align}
\end{assum}

One can show that Assumption~\ref{assum:sym_smooth} is equivalent to norm of Hessian being bounded by $\left (L_0+L_1\|\nabla f(y)\|_2 \right)$ provided the function is twice differentiable. We note that Assumption~\ref{assum:sym_smooth}, does not technically require the function to be differentiable twice. In the next section we will reference and also prove some properties of  $(L_0, L_1)$-smooth (convex) functions that will be later used to derive upper bounds.

\begin{remark} Assumption~\ref{assum:sym_smooth} follows from Assumption 1.2 from \cite{lobanov2025linearconvergencerateconvex} with the same $L_0$ and $L_1$. We also refer to \cite{chen_generalized_smooth_2023} for other propositions that are equivalent to $(L_0,L_1)$-smoothness as it is defined in this paper.
\end{remark}

\section{Properties of $(L_0, L_1)$-smooth convex functions}

We begin with a well-known equivalent definition of smooth convex class $\mathcal{F}_L$ \cite{Nesterov2013-mm}:
\begin{teo} [Equivalent definition for $\mathcal{F}_L$] \label{teo:smooth_convex_interpolation} Assumption~\ref{assum:convex} ($\mu = 0$) and Assumption~\ref{assum:smooth} are satisfied if and only if
\begin{align}
    f(y) \geqslant f(x) + \langle \nabla f(x), y - x\rangle + \frac{\|\nabla f(y)-\nabla f(x) \|_2^2}{2L}, \hspace{10pt} \forall x,y \in \R^d.
\end{align}
\end{teo}

The inequalities of Theorem~\ref{teo:smooth_convex_interpolation} are actually necessary and sufficient for interpolation \cite{taylor2017smooth} which can be used to obtain exact worst-case guarantees via PEP \cite{goujaud2023fundamental}. They also do not ruin a semidefinite problem because they involve distances and gradients as a quadratic form. There are similar inequalities in literature for symmetric $(L_0,L_1)$-smooth functions \cite{lobanov2025linearconvergencerateconvex}:

\begin{lem} \label{lem:symm_constraints} Let function $f$ satisfy Assumption~\ref{assum:sym_smooth} and Assumption~\ref{assum:convex} (convexity, $\mu=0$), then for $x,y \in \sR$ such that $\|y-x \|_2 \leq \frac{1}{L_1}$ we have:
\begin{align}
    \frac{\|\nabla f(x)-\nabla f(y)\|_2^2}{2(L_0+L_1\|\nabla f(y)\|_2)} \leqslant f(y) - f(x) - \left\langle \nabla f(x), y-x\right\rangle.
\end{align}
\end{lem}

Because Lemma~\ref{lem:symm_constraints} and other similar results (Lemma 2.2. \cite{gorbunov2025methods}) include bounded distance $\|y-x\|_2$, they are not convenient for our purposes. We shall derive our own proposition that works for arbitrarily spaced $x$ and $y$. First, we refer to the following properties of symmetric $(L_0, L_1)$ smooth functions:  

\begin{lem}[(Proposition~1 from \cite{chen_generalized_smooth_2023})] \label{lem:sym_equiv} Assumption~\ref{assum:sym_smooth} holds if and only if
\begin{align}
    \| \nabla f(y) - \nabla f(x) \|_2 \leqslant \left(L_0+L_1\|\nabla f(x)\|_2\right) \exp{(L_1\|x-y\|_2)}\|x-y\|_2, \hspace{10pt} \forall x,y \in \R^d.
\end{align}
Moreover Assumption~\ref{assum:convex} implies \cite{gorbunov2025methods} for all $x,y \in \R^d$
\begin{align}
    f(y) \leqslant f(x) + \left\langle \nabla f(x), y-x \right\rangle + \frac{L_0+L_1\|\nabla f(x)\|}{2}\exp{(L_1\|x-y\|_2)}\|x-y\|_2^2.
\end{align}
\end{lem}

With Lemma~\ref{lem:sym_equiv}, we can derive the required inequality.

\begin{lem} \label{lem:sym_conditions} If Assumption~\ref{assum:convex} and Assumption~\ref{assum:sym_smooth} hold then for any $x,y \in R^d$ the following inequality is satisfied:
\begin{align}
\label{eq:cond1}
    \frac{\|\nabla f(x)-\nabla f(y)\|_2^2}{2\left(L_0+L_1\|\nabla f(y)\|_2\right)\exp{\left(L_1\|x-y \|_2\right)}} \leqslant f(y) - f(x) - \left\langle \nabla f(x), y-x\right\rangle.
\end{align}
\end{lem}
\proof
This Lemma's proof is based on the proof of Lemma 2.2. from \cite{gorbunov2025methods}. First we introduce $\phi_x(y) := f(y) - \langle \nabla f(x), y\rangle$. By using Lemma~\ref{lem:sym_equiv} and doing the same calculations \cite{gorbunov2025methods}, we get $\forall u,y \in \R^d$:
\begin{align}
\label{eq:phi_cond}
    \phi_x(y) \leqslant \phi_x(u) + \langle \nabla\phi_x(u), y-u\rangle + \frac{L_0+L_1\|\nabla f(u)\|_2}{2} \exp{(L_1\|u-y\|_2)}\|u-y \|_2^2.
\end{align}
Next we introduce $\nu:= \exp{\left(-L_1\|x-u \|_2\right)}$ and derive the following inequality
\begin{align}
\nu\exp{\left(\frac{\nu L_1 \|\nabla\phi_x(u)\|_2}{L_0+L_1\|\nabla f(u)\|_2}\right)} & = 
\nu\exp{\left(\frac{\nu L_1 \|\nabla f(u) - \nabla f(x)\|_2}{L_0+L_1\|\nabla f(u)\|_2}\right)} \\ &\leqslant \nu \exp{\left(\frac{\nu L_1 \|x-u\|_2 \exp{(L_1 \|x-u\|_2)}(L_0+L_1\|\nabla f(u)\|_2)}{L_0+L_1\|\nabla f(u)\|_2}\right)} \\&=
\nu \exp{\left(\nu L_1 \|x-u\|_2 \exp{(L_1 \|x-u\|_2)}\right)} \\& = 1.
\end{align}
Using the bound above and Eq~\ref{eq:phi_cond} with $y = u - \frac{\nu}{L_0+L_1\|\nabla f(u)\|_2}\nabla\phi_x(u)$, we derive 
\begin{align}
\phi_x &\left( u - \frac{\nu}{L_0+L_1 \|\nabla f(u)\|_2}\nabla\phi_x(u) \right)  \\
&\leqslant \phi_x(u) - \nu \frac{\|\nabla\phi_x(u)\|_2^2}{\left(L_0+L_1 \|\nabla f(u)\|_2\right)} + \frac{\nu ^2\|\nabla\phi_x(u)\|_2^2}{2\left(L_0+L_1 \|\nabla f(u)\|_2\right)}\exp{\left(\frac{\nu L_1 \|\nabla\phi_x(u)\|_2}{L_0+L_1\|\nabla f(u)\|_2}\right)}\\
&\leqslant  \phi_x(u) - \nu \frac{\|\nabla\phi_x(u)\|_2^2}{\left(L_0+L_1 \|\nabla f(u)\|_2\right)} + \nu \frac{\|\nabla\phi_x(u)\|_2^2}{2\left(L_0+L_1 \|\nabla f(u)\|_2\right)}\\
&\leqslant \phi_x(u) - \nu \frac{\|\nabla\phi_x(u)\|_2^2}{2\left(L_0+L_1 \|\nabla f(u)\|_2\right)}.
\end{align}
Using the definition of $\phi_x(u)$ and the fact that $x$ is the minimum of $\phi_x(u)$, it follows that
\begin{align*}
    f(x) - \langle \nabla f(x),x\rangle \leqslant f(u) - \langle \nabla f(x), u \rangle -\frac{\nu \|\nabla f(x) - \nabla f(u)\|_2^2}{2\left(L_0+L_1\|\nabla f(u)\|_2\right)}, \hspace{4pt} x,u \in \R^d.
\end{align*}
which is equivalent to the statement of the Lemma.

\endproof

It's easy to see that when applied to every pair $(x,y)$, \eqref{eq:cond1} imply both $(L_0,L_1)$-smoothness and convexity. Therefore Lemma~\ref{lem:sym_conditions} can be seen as an equivalent definition of $(L_0,L_1)$-smooth convex functions. However, when applied to finite number of points, \eqref{eq:cond1} will (possibly) not be sufficient for interpolation by such function. This is acceptable for the purpose of this paper since our goal is to obtain better upper bounds and we are not aiming for exact worst-case bounds \cite{rubbens2023interpolation} in our approach. Besides, in order to make \eqref{eq:cond1} compatible with a semidefinite program, we will have relax the inequalities anyway. Next, we introduce standard Gradient Descent method (GD) and also $(L_0,L_1)-GD$ that will be analyzed in this paper:

\noindent
\begin{minipage}[t]{0.48\textwidth} 
\begin{algorithm}[H]
\caption{Gradient descent method.}
\label{alg:gd}
\begin{algorithmic}[1]
\STATE \textbf{Input:} Starting point $x_{0}$, number of steps $N$, stepsize $\eta$, $L > 0$
\FOR {$k = 1, \dots, N$}
    \STATE $x^{k+1} = x^{k} - \frac{\eta}{L} \nabla f(x^k)$
\ENDFOR
\STATE \textbf{Output:} $x^N$
\end{algorithmic}
\end{algorithm}
\end{minipage}
\hfill
\begin{minipage}[t]{0.48\textwidth} 
\begin{algorithm}[H]
\caption{$(L_0,L_1)$-GD.}
\label{alg:l0l1_gd}
\begin{algorithmic}[1]
\STATE \textbf{Input:} Starting point $x^{0}$, number of steps $N$, stepsize $\eta$, $L_0 > 0$, $L_1 \geqslant 0$
\FOR {$k = 1, \dots, N$}
    \STATE $x^{k+1} = x^{k} - \frac{\eta}{L_0+L_1\|\nabla f(x^k)\|_2} \nabla f(x^k)$
\ENDFOR
\STATE \textbf{Output:} $x^N$
\end{algorithmic}
\end{algorithm}
\end{minipage}
\noindent

It's known that as long as the gradient norm exceeds a certain threshold, $(L_0,L_1)$-GD exhibits linear convergence \cite{gorbunov2025methods}. From practical standpoint, it may be useful to know how algorithms behave under special conditions like gradient norm regime that we just mentioned. Thus, in this paper we will employ the novel approach to improve the upper bound for function error's contraction rate after one step of $(L_0,L_1)$-GD when $\|\nabla f(x^k)\|_2 > L_0/L_1$. As we suggested in the introduction to the paper, additional knowledge of system's dynamics can be used to derive better relaxations for inequalities characterizing the system.
In the case of $(L_0,L_1)$-GD, we can refer to the following property of gradient norm monotonicity:

\begin{lem} [Monotonicity of gradient norm \cite{gorbunov2025methods}] \label{lem:monotone_grad} Let Assumption~\ref{assum:convex} with $\mu = 0$ and Assumption~\ref{assum:sym_smooth} hold. Then for all $k \geqslant 0$ the iterates generated by $(L_0,L_1)$-GD with $\eta \leqslant \nu$, $\nu = \exp{(-\nu)}$ satisfy
\begin{align*}
    \|\nabla f(x^{k+1})\|_2 \leqslant \|\nabla f(x^k) \|_2.
\end{align*}
\end{lem}

With that, we can move to the essence of the paper.

\section{Construction of the PEP template for $(L_0, L_1)$-smooth convex functions} \label{sec:construction}

The idea of our proposed approach is to use existing PEP solutions as a templates to derive guarantees for systems that do not have an SDP-compatible characterization. To this end, we have to (first) prepare the template itself. One can notice that the inequalities from Lemma~\ref{lem:sym_conditions} possess a structure that is seemingly similar to necessary and sufficient interpolation conditions for smooth convex functions (Theorem~\ref{teo:smooth_convex_interpolation}). Let the constraints from Lemma~\ref{lem:sym_conditions} be imposed on finite pairs of points $(x^i,x^j)$. By letting $L = \max_i \left[ \ L_0 + L_1 \|\nabla f(x^i)\|_2\exp{\left(L_1\|x^i-x^j \|_2\right)}\right]$ and swapping the corresponding parts for $L$ in every constraint, we get a relaxed set of inequalities. Also, this new set will be the exactly as the set of interpolation constraints from Theorem~\ref{teo:smooth_convex_interpolation} applied to the same pairs of points. By redefining steps of $(L_0,L_1)$-GD, in a similar manner, we get a standard gradient descent method with constant step size $\eta^*$. The corresponding Performance estimation problem (for function value error as a performance metric) was accurately solved in \cite{taylor2017smooth} for step sizes $0 < \eta < 2$ and an arbitrary number of steps $N$. Unlike \cite{taylor2017smooth} that only bounds initial distance to minimum $\|x^0 - x_\ast\|_2 \leqslant R$, in this paper we shall introduce additional constraints setting initial function error $f(x^0)-f(x_*) = F$ and $\|\nabla f(x^0)\|_2 = G$. We shall also set $N = 1$ to reduce conservatism. \par

However, because a template based on minimization of a smooth convex function is too crude in the particular case of $(L_0, L_1)$-smoothness, so we defer it to the Appendix~\ref{appendix:smooth_template}. Here, instead of bounding $L_0 + L_1 \|\nabla f(x)\|_2$ parts with the same $L$ everywhere, we could get a better relaxation by using separate constant for every inequality. Thus, we write \eqref{eq:cond1} for all pairs of $x^k$, $x^{k+1}$ and $x_*$:
\begin{align}
    &f(x_*) \geqslant f(x^k) + \langle \nabla f(x^k),x_{*}-x^k\rangle + \frac{\|\nabla f(x^k)-\nabla f(x_*) \|_2^2}{2 l_{0*}}, \label{eq:conditions} \\
    &f(x_*) \geqslant f(x^{k+1}) + \langle \nabla f(x^{k+1}),x_{*}-x^{k+1}\rangle + \frac{\|\nabla f(x^{k+1})-\nabla f(x_*) \|_2^2}{2 l_{1*}}, \notag \\
    &f(x^k) \geqslant f(x_*) + \langle \nabla f(x_*),x^k-x_*\rangle + \frac{\|\nabla f(x_*)-\nabla f(x^k) \|_2^2}{2 l_{*0}}, \notag \\
    &f(x^k) \geqslant f(x^{k+1}) + \langle \nabla f(x^{k+1}),x^k-x^{k+1}\rangle + \frac{\|\nabla f(x^{k+1})-\nabla f(x^k) \|_2^2}{2 l_{10}}, \notag \\
    &f(x^{k+1}) \geqslant f(x_*) + \langle \nabla f(x_*),x^{k+1}-x_*\rangle + \frac{\|\nabla f(x_*)-\nabla f(x^{k+1}) \|_2^2}{2 l_{*1}}, \notag \\
    &f(x^{k+1}) \geqslant f(x^k) + \langle \nabla f(x^k),x^{k+1}-x^k\rangle + \frac{\|\nabla f(x^k)-\nabla f(x^{k+1}) \|_2^2}{2 l_{01}}. \notag
\end{align}
where, assuming we make one step of $(L_0,L_1)$-GD: $x^{k+1} = x^{k} - \frac{\eta}{L_0+L_1\|\nabla f(x^k)\|_2} \nabla f(x^k)$ and $\|x^k-x_*\|_2 \leqslant R$, $\|\nabla f(x^k)\|_2 = G$, we have
\begin{align*}
    &l_{0*} := \left(L_0+L_1\|\nabla f(x_*)\|_2\right)\exp{\left(L_1\|x^k-x_* \|_2\right)} = L_0\exp{(L_1R)}, \\
    &l_{*0}  := \left(L_0+L_1\|\nabla f(x^k)\|_2\right)\exp{\left(L_1\|x^k-x_* \|_2\right)} = \left(L_0+L_1G\right)\exp{(L_1R)}, \\
    &l_{1*} := \left(L_0+L_1\|\nabla f(x_*)\|_2\right)\exp{\left(L_1\|x^{k+1}-x_* \|_2\right)} \leqslant L_0\exp{(L_1R)}, \\
    &l_{*1} := \left(L_0+L_1\|\nabla f(x^{k+1})\|_2\right)\exp{\left(L_1\|x^{k+1}-x_* \|_2\right)} \leqslant \left(L_0+L_1G\right)\exp{(L_1R)}, \\
    &l_{10} := \left(L_0+L_1\|\nabla f(x^k)\|_2\right)\exp{\left(L_1\|x^k-x^{k+1} \|_2\right)} = \left(L_0+L_1G\right)\exp{\left(\frac{\eta GL_1}{L_0+L_1G}\right)}, \\
    &l_{01} := \left(L_0+L_1\|\nabla f(x^{k+1})\|_2\right)\exp{\left(L_1\|x^k-x^{k+1} \|_2\right)} \leqslant \left(L_0+L_1G\right)\exp{\left(\frac{\eta GL_1}{L_0+L_1G}\right)}.
\end{align*}

Given the constants listed above, one would need to formulate and solve a PEP with 3 "smoothness" parameters: $l_{0*}$($= l_{1*}$), $l_{*0}$($= l_{*1}$) and $l_{10}$($= l_{01}$). There is also step size constant and three initial conditions parameters which would lead to seven parameters in total. 
We can reduce the number of parameters of the problem. First, we choose the largest admissible step size $\eta = \nu$ that guarantees that gradient's norm is non-increasing (Lemma~\ref{lem:monotone_grad}):
\begin{align*}
    l_{10}, l_{01} &\leqslant  \left(L_0+L_1G\right)\exp{\left(\frac{\nu GL_1}{L_0+L_1G}\right)} \\ &\leqslant \frac{\left(L_0+L_1G\right)}{\nu}\nu\exp{\left(\nu\right)} = \frac{\left(L_0+L_1G\right)}{\nu}.
\end{align*}
So we can let 
\begin{align*}
    &l_0 :=  L_0\exp{(L_1R)},\\
    &l_1 := \frac
    {\left(L_0+L_1G\right)}{\nu}, \\
    &l_2 := \left(L_0+L_1G\right)\exp{(L_1R)}.
\end{align*}
and get the following relaxation of the initial system:
\begin{align}
    &f(x_*) \geqslant f(x^k) + \langle \nabla f(x^k),x_{*}-x^k\rangle + \frac{\|\nabla f(x^k)-\nabla f(x_*) \|_2^2}{2 l_{0}}, \label{eq:conditions_1} \\
    &f(x_*) \geqslant f(x^{k+1}) + \langle \nabla f(x^{k+1}),x_{*}-x^{k+1}\rangle + \frac{\|\nabla f(x^{k+1})-\nabla f(x_*) \|_2^2}{2 l_{0}}, \notag \\
    &f(x^k) \geqslant f(x_*) + \langle \nabla f(x_*),x^k-x_*\rangle + \frac{\|\nabla f(x_*)-\nabla f(x^k) \|_2^2}{2 l_{2}}, \notag \\
    &f(x^k) \geqslant f(x^{k+1}) + \langle \nabla f(x^{k+1}),x^k-x^{k+1}\rangle + \frac{\|\nabla f(x^{k+1})-\nabla f(x^k) \|_2^2}{2 l_{1}}, \notag \\
    &f(x^{k+1}) \geqslant f(x_*) + \langle \nabla f(x_*),x^{k+1}-x_*\rangle + \frac{\|\nabla f(x_*)-\nabla f(x^{k+1}) \|_2^2}{2 l_{2}}, \notag \\
    &f(x^{k+1}) \geqslant f(x^k) + \langle \nabla f(x^k),x^{k+1}-x^k\rangle + \frac{\|\nabla f(x^k)-\nabla f(x^{k+1}) \|_2^2}{2 l_{1}}, \notag
\end{align}
and the GD step: $x^{k+1} = x^{k} - \nabla f(x^k)/l_1$. Now, we introduce the semidefinite problem:

%More general PEP problem that can be used as framework to study (L_0,L_1)-smooth and other function classes
\begin{align*} \label{problem:G_PEP}\tag{G-PEP}
    \und{\sup}{\{x^i,\nabla f(x^i), f(x^i)\}_{i\in I}} &f(x^{k+1}) - f(x_*) \quad \text{such that}  \\
    &I = \{k, k+1,*\},\\
    & \{x^i,\nabla f(x^i), f(x^i)\}_{i\in I} \quad \text{satisfy \eqref{eq:conditions_1}} , \\
    &x^{k+1} = x^k - \frac{\nabla f(x^i)}{l_1},\\
    &\{x_*,\nabla f(x_*), f(x_*)\} = \{\textbf{0}, \textbf{0}, 0\},\\
    & \|x^k-x_*\|_2 \leqslant R, \\
    & \|\nabla f(x^k)\|_2 = G, \\
    & f(x^k)-f_* = F.
\end{align*}

We can establish some necessary conditions for feasibility of \ref{problem:G_PEP}:

%Lemma: necessary conditions for feasibility of G-PEP
\begin{lem} \label{lem:feasibility_ness_gen1} For \ref{problem:G_PEP} to be feasible, it is necessary that
\begin{align*}
    &G \leqslant 2 \frac{l_0l_1}{l_0+l_1}R, \\
    &\frac{G^2}{2l_1} \leqslant F \leqslant RG - \frac{G^2}{2l_0}.
\end{align*}
\end{lem}
\proof

Consider the following constraints from~\ref{problem:G_PEP}:
\begin{align*}
    &\lambda_{s0}: f(x_*) \geqslant f(x^k) + \langle \nabla f(x^k), x_{*}-x^k\rangle + \frac{\|\nabla f(x^k)-\nabla f(x_*) \|_2^2}{2l_0}, \\
    &\lambda_{01}: f(x^{k}) \geqslant f(x^{k+1}) + \langle \nabla f(x^{k+1}), x^{k}-x^{k+1} \rangle + \frac{\|\nabla f(x^{k+1})-\nabla f(x^k) \|_2^2}{2l_1}, \\
    & \lambda_R:R^2 - \|x^{k}-x_{*} \|_2^2 \geqslant 0, \\
    &\lambda_G: G^2 - \|\nabla f(x^k) \|_2^2 = 0. \\
\end{align*}
First we prove the condition for $G$. By adding up the constraints with the following multipliers
\begin{align*}
    &0 \leqslant \lambda_{s0} = \frac{4l_1}{1+l_1/l_0}, \\
    &0 \leqslant \lambda_{01} = \lambda_{s0},\\
    &\lambda_{R} = \lambda_{s0} \frac{l_1}{1+l_1/l_0}.
\end{align*}
and substituting $x^{k+1} = x^{k} - \nabla f(x^k)/l_1$, we get 
\begin{align*}
    G^2 = \|\nabla f(x^k)\|_2^2 \leqslant  \frac{4l_1^2R^2}{\left(1+l_1/l_0\right)^2} - \lambda_{s0}\left(f(x^{k+1}) - f(x_*)\right) - p^\top (S\otimes I_n) p ,
\end{align*}
where $p:=((x^k-x_{*})^\top, \left(\nabla f(x^k)-\nabla f(x_*))^\top,(\nabla f(x^{k+1})-\nabla f(x_*))^\top \right)^\top$ and 
\begin{align*}
    &S:= \begin{bmatrix*}
	\frac{4l_1^2}{\left(1+l_1/l_0\right)^2}  & -\frac{2l_1}{1+l_1/l_0} & 0 \\
	-\frac{2l_1}{1+l_1/l_0} & 1 & 0 \\
	0  & 0 & \frac{2l_1}{1+l_1/l_0}
    \end{bmatrix*}
\\
\end{align*}
Next, we prove the conditions for $F$. Using
\begin{align*}
    &0 \leqslant\lambda_{01} = 1, \\
    &\lambda_{G} = -\frac{1}{2l_1},
\end{align*}
we get the lower bound
\begin{align*}
    F = f(x^k) - f(x_*) &\geqslant \frac{G^2}{2l_1} + \left(f(x^{k+1}) - f(x_*)\right) + \frac{\|\nabla f(x^{k+1})\|_2^2}{2l_1} \\
    &\geqslant \frac{G^2}{2l_1}.
\end{align*}
By adding up the constraints with the following multipliers
\begin{align*}
    &0 \leqslant\lambda_{s0} = 1, \\
    &0 \leqslant \lambda_{R} = \frac{G}{2R},\\
    &\lambda_{G} = \frac{R}{2G} - \frac{1}{2l_0},
\end{align*}
we get the upper bound
\begin{align*}
    F =f(x^k) - f(x_*) \leqslant GR - \frac{G^2}{2l_0} - p^\top (S\otimes I_n) p,
\end{align*}
where $p:=((x^k-x_{*})^\top, \left(\nabla f(x^k)-\nabla f(x_*))^\top \right)^\top$ and 
\begin{align*}
    &S:= \begin{bmatrix*}
	\lambda_R  & -1/2  \\
	-1/2 & \left(\frac{1}{2l_0}+\lambda_G\right) 
    \end{bmatrix*}.
\\
\end{align*}

\endproof

We verified numerically (Appendix~\ref{appendix:pep_experiments}) that the conditions from Lemma~\ref{lem:feasibility_ness_gen1} are actually sufficient for feasibility of~\ref{problem:G_PEP} with one exception. The bound $F \leqslant RG - \frac{G^2}{2l_0}$ is not sufficient for problem's feasibility when $l_0>l_1$ and the bound could be strengthened. The closed form expression for necessary and sufficient bound appears to be extremely complicated based on preliminary analysis, so we do not provide it in this paper. Nevertheless, we managed to verify that relative difference between the true bound $F^{true}$ and our estimate $F^{est} = RG - \frac{G^2}{2l_0}$ is not big:
\begin{align*}
    \frac{F^{est} - F^{true}}{F^{est}} \leqslant 0,065.
\end{align*}
Therefore, we might as well be satisfied with the estimate from Lemma~\ref{lem:feasibility_ness_gen1}. 

The problem~\ref{problem:G_PEP} is essentially a single-step PEP for GD with iterates constrained by \ref{eq:conditions_1}. Let \text{G-Val} be the value of the problem~\ref{problem:G_PEP} In the next Theorem we derive an upper bound for \text{G-Val}.

%Theorem: upper bound for problem G-PEP (l0 >= l1)
\begin{teo} \label{teo:general_FG_res} Consider problem~\ref{problem:G_PEP} with $l_2 \geqslant l_0 \geqslant l_1$. Suppose $R,G,F > 0$ satisfy Lemma~\ref{lem:feasibility_ness_gen1}. If the problem is feasible then for it's value, \text{G-Val}, we state that 
\begin{align}
    f(x^{k+1})-f(x_*) \leqslant \text{G-Val} \leqslant \left(F -\frac{G^2}{2l_1} \right)\left(1-\frac{F -\frac{G^2}{2l_1}}{4\left[\frac{l_1R^2}{2}+\frac{G^2}{8l_1}\left(1-(l_1/l_0)^2\right)-\frac{F}{2}\left(1+l_1/l_0\right)\right]}\right).
    \label{eq:g_val_bound}
\end{align}
if
\begin{align}
    &G \leqslant 2 \frac{l_0l_1}{l_0+l_1}R, \\
    &\frac{G^2}{2l_1} \leqslant F < \min \left[\alpha,\beta\right]  ,
    \label{eq:f_cond_l0l1}
\end{align}
where
\begin{align*}
    &\alpha := \frac{1}{3}\left[\frac{2l_1R^2}{1+l_1/l_0}+\frac{G^2}{2l_1}\left(2-\frac{l_1}{l_0}\right) \right], \\
    &\beta := RG - \frac{G^2}{2l_0}.
\end{align*}
otherwise if
\begin{align*}
    &G \leqslant 2 \frac{l_0l_1}{l_0+l_1}R, \\
    &\frac{1}{3}\left[\frac{2l_1R^2}{1+l_1/l_0}+\frac{G^2}{2l_1}\left(2-\frac{l_1}{l_0}\right) \right]\leqslant F \leqslant RG - \frac{G^2}{2l_0}.
\end{align*}
then
\begin{align}
    f(x^{k+1}) - f(x_*) \leqslant \text{G-Val} \leqslant\frac{2l_1R^2}{\left(1+l_1/l_0\right)^2} + \frac{G^2}{l_1}\frac{1-l_1/l_0}{1+l_1/l_0} + \left(1-\frac{4}{1+l_1/l_0}\right)F.
\label{eq:g_val_bound_1}
\end{align}
\end{teo}

\proof
Consider the following constraints from~\ref{problem:G_PEP}:
\begin{align*}
    &\lambda_{s0}: f(x_*) \geqslant f(x^k) + \langle \nabla f(x^k), x_{*}-x^k\rangle + \frac{\|\nabla f(x^k)-\nabla f(x_*) \|_2^2}{2l_0}, \\
    &\lambda_{s1}: f(x_*) \geqslant f(x^{k+1}) + \langle \nabla f(x^{k+1}),x_{*}-x^{k+1}\rangle + \frac{\|\nabla f(x^{k+1})-\nabla f(x_*) \|_2^2}{2l_0}, \\
    &\lambda_{01}: f(x^{k}) \geqslant f(x^{k+1}) + \langle \nabla f(x^{k+1}), x^{k}-x^{k+1} \rangle + \frac{\|\nabla f(x^{k+1})-\nabla f(x^k) \|_2^2}{2l_1}, \\
    &\lambda_{10}: f(x^{k+1}) \geqslant f(x^{k}) + \langle \nabla f(x^{k}), x^{k+1}-x^{k} \rangle + \frac{\|\nabla f(x^{k+1})-\nabla f(x^k) \|_2^2}{2l_1}, \\
    & \lambda_R:R^2 - \|x^{k}-x_{*} \|_2^2 \geqslant 0, \\
    &\lambda_G: G^2 - \|\nabla f(x^k)\|_2^2 = 0, \\
    &\lambda_F: F-\left(f(x^k)-f(x_*)\right) = 0.
\end{align*}
We add them up with following multipliers:
\begin{align*}
    &0 \leqslant \lambda_{s1} = -\frac{1}{2} \frac{\frac{G^2}{2l_1}-F}{\frac{l_1R^2}{2}+\frac{G^2}{8l_1}\left(1-(l_1/l_0)^2\right)-\frac{F}{2}\left(1+l_1/l_0\right)},\\
    &0 \leqslant\lambda_{s0} = \frac{\lambda_{s1}^2}{2}\left(1+\frac{l_1}{l_0}\right), \\
    &0 \leqslant \lambda_{10} = \frac{\lambda_{s1}}{2}\left(1-\frac{l_1}{l_0}\right), \\
    &0 \leqslant \lambda_{01} = 1 - \lambda_{s1} + \lambda_{10}, \\
    &\lambda_{F} = 1 - \lambda_{s1} - \lambda_{s0}, \\
    &0 \leqslant \lambda_{R} = \frac{l_1}{2}\lambda_{s1}^2, \\
    &\lambda_{G} = \frac{\lambda_{s1}-1}{2l_1}+\frac{\lambda_{s1}^2}{8l_1}\left(1-\left(\frac{l_1}{l_0}\right)^2 \right)
\end{align*}
while also substituting $x^{k+1} = x^{k} - \nabla f(x^k)/l_1$:
\begin{align}
    f(x^{k+1}) - f(x_*) \leqslant \left(F -\frac{G^2}{2l_1} \right)\left(1-\frac{F -\frac{G^2}{2l_1}}{4\left[\frac{l_1R^2}{2}+\frac{G^2}{8l_1}\left(1-(l_1/l_0)^2\right)-\frac{F}{2}\left(1+l_1/l_0\right)\right]}\right)- p^\top (S\otimes I_n) p,
\label{eq:gen_FG_res_2}
\end{align}
where $p:=((x^k-x_{*})^\top, \left(\nabla f(x^k)-\nabla f(x_*))^\top,(\nabla f(x^{k+1})-\nabla f(x_*))^\top \right)^\top$ and 
\begin{align*}
    &S:= \frac{1}{2} \begin{bmatrix*}
	l_1\lambda_{s1}^2  & -\frac{\lambda_{s1}^2}{2}\left(1+l_1/l_0\right) & -\lambda_{s1} \\
	-\frac{\lambda_{s1}^2}{2}\left(1+l_1/l_0\right) & \frac{\lambda_{s1}^2}{4l_1}\left(1+l_1/l_0\right)^2 & \frac{\lambda_{s1}-\lambda_{10}}{l_1} \\
	-\lambda_{s1}  & \frac{\lambda_{s1}-\lambda_{10}}{l_1} & \frac{1+2\lambda_{10}}{l_1}-\frac{\lambda_{s1}}{l_1}\left(1-l_1/l_0\right)
    \end{bmatrix*}
\\
\end{align*}
is a positive semidefinite matrix. We have to ensure that multipliers for the inequalities are non-negative. From $\lambda_{s1} \geqslant 0$:
\begin{align*}
    &F \geqslant \frac{G^2}{2l_1}, \\
    &F \leqslant \frac{l_0l_1R^2}{l_0+l_1} + \frac{G^2}{4l_1} \left(1-\frac{l_1}{l_0}\right),
\end{align*}
it follows that
\begin{align}
    G \leqslant 2 \frac{l_0l_1}{l_0+l_1}R.
\label{eq:G_restr}
\end{align}
From $\lambda_{01} \geqslant 0$ we get
\begin{align}
    F \leqslant \frac{2}{3(1+l_1/l_0)}\left[l_1R^2+\frac{G^2}{4l_1}\left(2+\frac{l_1}{l_0} -\left(\frac{l_1}{l_0}\right)^2\right)\right].
\label{eq:F_restr}
\end{align}
Because for all $G$ that satisfy \eqref{eq:G_restr}:
\begin{align*}
    \frac{l_0l_1R^2}{l_0+l_1} + \frac{G^2}{4l_1} \left(1-\frac{l_1}{l_0}\right) \geqslant \frac{2}{3(1+l_1/l_0)}\left[l_1R^2+\frac{G^2}{4l_1}\left(2+\frac{l_1}{l_0} -\left(\frac{l_1}{l_0}\right)^2\right)\right].
\end{align*}
we only need the second inequality for $F$, \eqref{eq:F_restr} which can be rewritten as
\begin{align}
    F \leqslant \frac{1}{3}\left[\frac{2l_1R^2}{1+l_1/l_0}+\frac{G^2}{2l_1}\left(2-\frac{l_1}{l_0}\right) \right]. 
\label{eq:F_restr_1}
\end{align}
Now consider a different set of multipliers:
\begin{align*}
    &0 \leqslant \lambda_{s1} = \frac{2}{1+l_1/l_0},\\
    &0 \leqslant\lambda_{s0} = \frac{2}{1+l_1/l_0}, \\
    &0 \leqslant \lambda_{10} = \frac{1-l_1/l_0}{1+l_1/l_0}, \\
    &\lambda_{F} = 1 - \frac{4}{1+l_1/l_0}, \\
    &0 \leqslant \lambda_{R} = \frac{2l_1}{\left(1+l_1/l_0\right)^2}, \\
    &\lambda_{G} = \frac{1-l_1/l_0}{l_1\left(1+l_1/l_0\right)}.
\end{align*}
The same way we get
\begin{align}
    f(x^{k+1}) - f(x_*) \leqslant \frac{2l_1R^2}{\left(1+l_1/l_0\right)^2} + \frac{G^2}{l_1}\frac{1-l_1/l_0}{1+l_1/l_0} + \left(1-\frac{4}{1+l_1/l_0}\right)F.
\label{eq:gen_FG_res_3}
\end{align}
The bound above is valid for all $F$, including 
\begin{align*}
    F > \frac{1}{3}\left[\frac{2l_1R^2}{1+l_1/l_0}+\frac{G^2}{2l_1}\left(2-\frac{l_1}{l_0}\right) \right],
\end{align*}
however, for $F$ satisfying \eqref{eq:F_restr_1}, the first bound \eqref{eq:gen_FG_res_2} is smaller.

\endproof

Now we shall consider the case $l_1 > l_0$.

%Theorem: upper bound for problem G-PEP (l0 < l1)
\begin{teo} \label{teo:general_FG_res_1} Consider problem~\ref{problem:G_PEP} with $l_2 \geqslant l_0 < l_1$. Suppose $R,G,F > 0$ satisfy Lemma~\ref{lem:feasibility_ness_gen1}. If the problem is feasible then for it's value, \text{G-Val}, we state that if
\begin{align}
    &G \leqslant 2 \frac{l_0l_1}{l_0+l_1}R, \\
    &\frac{G^2}{2l_1} \leqslant F < \min \left[d_1, \ d_2\right],  
    \label{eq:f_cond_l0l1_2}
\end{align}
where
\begin{align*}
    &d_1 := \frac{l_1/l_0+1}{\left(l_1/l_0\right)^2+l_1/l_0+1}\left[\frac{l_1R^2}{2} + \frac{G^2}{2l_1}\frac{1}{l_1/l_0+1}\right], \\
    &d_2 := RG - \frac{G^2}{2l_0}
\end{align*}
then
\begin{align}
    f(x^{k+1})-f(x_*) \leqslant \text{G-Val} \leqslant \sqrt{\frac{\alpha}{\beta}}\left(\frac{\sqrt{\alpha}-\sqrt{\beta}}{l_1/l_0-1}\right)^2 + \frac{\gamma}{l_1/l_0-1}\left(\frac{l_1}{l_0} - \sqrt{\frac{\alpha}{\beta}}\right),
    \label{eq:g_val_bound_2}
\end{align}
where
\begin{align*}
    &\alpha := \frac{l_1R^2}{2}-\frac{G^2}{2l_1}\left(\frac{l_1}{l_0}-1\right) - F, \\
    &\beta := \frac{l_1R^2}{2} - \frac{l_1}{l_0}F, \\
    &\gamma = F - \frac{G^2}{2l_1}.
\end{align*}
Otherwise, if
\begin{align*}
    &G \leqslant 2 \frac{l_0l_1}{l_0+l_1}R, \\
    &\frac{l_1/l_0+1}{\left(l_1/l_0\right)^2+l_1/l_0+1}\left[\frac{l_1R^2}{2} + \frac{G^2}{2l_1}\frac{1}{l_1/l_0+1}\right]< F \leqslant RG - \frac{G^2}{2l_0}
\end{align*}
then
\begin{align}
    f(x^{k+1}) - f(x_*) \leqslant \text{G-Val} \leqslant \frac{l_0R^2}{2} - \frac{G^2}{2l_1}\left(1-l_0/l_1\right) -\frac{l_0}{l_1}F
\label{eq:g_val_bound_3}
\end{align}
\end{teo}

\proof

Consider the following constraints from~\ref{problem:G_PEP}:
\begin{align*}
    &\lambda_{s0}: f(x_*) \geqslant f(x^k) + \langle \nabla f(x^k), x_{*}-x^k\rangle + \frac{\|\nabla f(x^k)-\nabla f(x_*) \|_2^2}{2l_0}, \\
    &\lambda_{s1}: f(x_*) \geqslant f(x^{k+1}) + \langle \nabla f(x^{k+1}),x_{*}-x^{k+1}\rangle + \frac{\|\nabla f(x^{k+1})-\nabla f(x_*) \|_2^2}{2l_0}, \\
    &\lambda_{01}: f(x^{k}) \geqslant f(x^{k+1}) + \langle \nabla f(x^{k+1}), x^{k}-x^{k+1} \rangle + \frac{\|\nabla f(x^{k+1})-\nabla f(x^k) \|_2^2}{2l_1}, \\
    & \lambda_R:R^2 - \|x^{k}-x_{*} \|_2^2 \geqslant 0, \\
    &\lambda_G: G^2 - \|\nabla f(x^k)\|_2^2 = 0,\\
    &\lambda_F: F-\left(f(x^k)-f(x_*)\right) = 0.
\end{align*}

We add them up with following multipliers:
\begin{align*}
    &0 \leqslant \lambda_{s1} = \frac{1}{l_1/l_0-1}\left[\sqrt{\frac{\frac{l_1R^2}{2}-\frac{G^2}{2l_1}\left(\frac{l_1}{l_0}-1\right) - F}{\frac{l_1R^2}{2} - \frac{l_1}{l_0}F}}-1\right],\\
    &0 \leqslant\lambda_{s0} = \frac{\lambda_{s1}^2}{1+\lambda_{s1}\left(l_1/l_0-1\right)}, \\
    &0 \leqslant \lambda_{01} = 1 - \lambda_{s1}, \\
    &\lambda_{F} = 1 - \lambda_{s1} - \lambda_{s0}, \\
    &0 \leqslant \lambda_{R} = \frac{l_1}{2}\lambda_{s0}, \\
    &\lambda_{G} = \frac{\lambda_{s1}-1}{2l_1} + \frac{\lambda_{s0}}{2l_1}\left[1-\frac{l_1}{l_0}\right]
\end{align*}
and substituting $x^{k+1} = x^{k} - \nabla f(x^k)/l_1$:
\begin{align*}
    f(x^{k+1}) - f(x_*) \leqslant \sqrt{\frac{\alpha}{\beta}}\left(\frac{\sqrt{\alpha}-\sqrt{\beta}}{l_1/l_0-1}\right)^2 + \frac{\gamma}{l_1/l_0-1}\left(\frac{l_1}{l_0} - \sqrt{\frac{\alpha}{\beta}}\right) - p^\top (S\otimes I_n) p,
\end{align*}
where 
\begin{align*}
    &\alpha := \frac{l_1R^2}{2}-\frac{G^2}{2l_1}\left(\frac{l_1}{l_0}-1\right) - F, \\
    &\beta := \frac{l_1R^2}{2} - \frac{l_1}{l_0}F, \\
    &\gamma = F - \frac{G^2}{2l_1}, \\
    &p:=\left((x^k-x_{*})^\top, \left(\nabla f(x^k)-\nabla f(x_*))^\top,(\nabla f(x^{k+1})-\nabla f(x_*)\right)^\top \right)^\top
\end{align*}
and 
\begin{align*}
    &S:= \begin{bmatrix*}
    \frac{l_1\lambda_{s0}}{2} & -\frac{\lambda_{s0}}{2} & -\frac{\lambda_{s1}}{2}\\
    -\frac{\lambda_{s0}}{2} & \frac{\lambda_{s0}}{2l_1} & \frac{\lambda_{s1}}{2l_1} \\
    -\frac{\lambda_{s1}}{2} & \frac{\lambda_{s1}}{2l_1} & \frac{\lambda_{s1}}{2l_0} + \frac{\lambda_{01}}{2l_1}
    \end{bmatrix*}
\\
\end{align*}
is a positive semidefinite matrix. We have to ensure that multipliers for the inequalities are non-negative. From $\lambda_{s1} \geqslant 0$ we have
\begin{align*}
    F \geqslant \frac{G^2}{2l_1}.
\end{align*}
Also, we have to make sure that the root in $\lambda_{s1}$ exists:
\begin{align*}
    F \leqslant \frac{l_0R^2}{2}.
\end{align*}
Both of the inequalities above hold automatically since we demanded that parameters satisfy Lemma~\ref{lem:feasibility_ness_gen1}. From $\lambda_{01} \geqslant 0$ we have 
\begin{align}
    F \leqslant \frac{l_1/l_0+1}{\left(l_1/l_0\right)^2+l_1/l_0+1}\left[\frac{l_1R^2}{2} + \frac{G^2}{2l_1}\frac{1}{l_1/l_0+1}\right].
    \label{eq:F_restr_2}
\end{align}
Now consider a different set of multipliers:
\begin{align*}
    &0 \leqslant \lambda_{s1} = 1,\\
    &0 \leqslant\lambda_{s0} = \frac{l_0}{l_1}, \\
    &0 \leqslant \lambda_{01} = 0, \\
    &\lambda_{F} = -\frac{l_0}{l_1}, \\
    &0 \leqslant \lambda_{R} = \frac{l_0}{2}, \\
    &\lambda_{G} = \frac{1}{2l_1}\left(\frac{l_0}{l_1}-1\right).
\end{align*}
Similarly, we obtain
\begin{align*}
    f(x^{k+1}) - f(x_*) \leqslant \frac{l_0R^2}{2} - \frac{G^2}{2l_1}\left(1-\frac{l_0}{l_1}\right)-\frac{l_0}{l_1}F.
\end{align*}
The bound above is valid for all $F$, including 
\begin{align*}
    F > \frac{l_1/l_0+1}{\left(l_1/l_0\right)^2+l_1/l_0+1}\left[\frac{l_1R^2}{2} + \frac{G^2}{2l_1}\frac{1}{l_1/l_0+1}\right].
\end{align*}
However, for $F$ satisfying \eqref{eq:F_restr_2}, the first bound \eqref{eq:g_val_bound_2} is smaller.

\endproof

Due to our deliberate choice of considering only one iteration in \ref{problem:G_PEP}, the template solutions should not be introducing too much conservatism into the resulting guarantees. The chosen parameterization is also versatile, so there is possibility of applying the proposed template to other definitions of generalized smoothness. For instance, it could be used for analyzing gradient descent variations in $(\rho, L_0, L_1)$-smoothness setting \cite{tyurin2025nearoptimalconvergenceacceleratedgradient} that is more general then $(L_0,L_1)$-smoothness considered in this paper. With that said, there is still room for improvement when it comes to versatility.

\section{Worst case contraction rate after one iteration}

%In this section we derive upper bounds for worst-case contraction rate of functiona 

%We begin with smooth convex functions:

In this section we find upper bound for contraction factor for one step of $(L_0,L_1)$-GD in the setting of symmetric $(L_0,L_1)$-smooth convex functions. First, a lemma for the case of $l_1 \leqslant l_0$

%Lemma for contraction factor upper bound: (L_0,L_1)-smooth convex functions (l0 >= l1)
\begin{lem} \label{lem:smooth_contraction_gen_1} Suppose the conditions of Theorem~\ref{teo:general_FG_res} are satisfied. For set initial gradient norm $G$ (and bounded distance to the optimum $R$), the following upper bound holds:
\begin{align}
    c = \frac{f(x^{k+1})-f(x_*)}{f(x^{k})-f(x_*)} &= \frac{f(x^{k+1})-f(x_*)}{F} \notag \\
    &\leqslant 1 - \frac{\theta}{\phi_*} - \frac{(\phi_*-\theta)^2}{\phi_*\left(4+\theta\left(1-(l_1/l_0)^2\right)-2\phi_*\left(1+l_1/l_0\right)\right)}. 
\label{eq:gen_contr_bound_1}
\end{align}
where 
\begin{align*}
    \theta := \left(\frac{G}{l_1R}\right)^2 \leqslant \frac{4}{\left(1+l_1/l_0\right)^2}
\end{align*}
and, depending on $\theta$, $\phi_*$ is equal to 
\begin{align*}
    &\textit{if} \quad  \theta \neq \frac{4}{\left(1+l_1/l_0\right)\left(7+5l_1/l_0\right)} \quad \textit{then} \\\\
    &\phi_* = \\\\ &\frac{\frac{\theta^2 l_1}{2l_0}\left(1-\frac{l_1(1+l_1/l_0)}{2l_0}\right) + \theta\left(\theta+2+2l_1/l_0\right)- \left(1-\theta (1+l_1/l_0)^2/4\right)\sqrt{\theta\left[\theta+4+\theta\left(1-(l_1/l_0)^2\right)\right]}}{\frac{\theta}{4}\left(1+l_1/l_0\right)\left(5l_1/l_0+7\right)-1} \\\\
    &\textit{if} \quad \theta = \frac{4}{\left(1+l_1/l_0\right)\left(7+5l_1/l_0\right)} \quad \textit{then} \quad \phi_* = \frac{4\left(2+l_1/l_0\right)}{\left(1+l_1/l_0\right)\left(7+5l_1/l_0\right)} \quad 
\end{align*}
\end{lem}
\proof
Consider the upper bound for $f(x^{k+1}) -f(x_*)$ from the statement of Theorem~\ref{teo:general_FG_res}. First assume $F > \frac{1}{3}\left[\frac{2l_1R^2}{1+l_1/l_0}+\frac{G^2}{2l_1}\left(2-\frac{l_1}{l_0}\right) \right]$. In this case which we call the second mode, the contraction expression is
\begin{align}
    \frac{f(x^{k+1}) - f(x_*)}{F} \leqslant \frac{1}{F}\left[\frac{2l_1R^2}{\left(1+l_1/l_0\right)^2} + \frac{G^2}{l_1}\frac{1-l_1/l_0}{1+l_1/l_0}\right] + 1-\frac{4}{1+l_1/l_0}.
\end{align}
From the expression above, the bound for contraction is a decreasing function of $F$ which attains maximum at $F = \frac{1}{3}\left[\frac{2l_1R^2}{1+l_1/l_0}+\frac{G^2}{2l_1}\left(2-\frac{l_1}{l_0}\right) \right]$. Therefore we only need to consider $F \leqslant \frac{1}{3}\left[\frac{2l_1R^2}{1+l_1/l_0}+\frac{G^2}{2l_1}\left(2-\frac{l_1}{l_0}\right) \right]$:
\begin{align}
    \frac{f(x^{k+1}) -f(x_*)}{F} \leqslant \left(1 -\frac{G^2}{2l_1F} \right)\left(1-\frac{F -\frac{G^2}{2l_1}}{4\left[\frac{l_1R^2}{2}+\frac{G^2}{8l_1}\left(1-(l_1/l_0)^2\right)-\frac{F}{2}\left(1+l_1/l_0\right)\right]}\right). \notag
\end{align}
We denote
\begin{align*}
    &\theta \colon= \left(\frac{G}{l_1R}\right)^2, \ \phi \colon= \frac{2F}{l_1R^2}.
\end{align*}
And rewrite right side expression as a function of $\phi$:
\begin{align*}
    c(\phi) := 1 - \frac{\theta}{\phi} - \frac{(\phi-\theta)^2}{\phi\left(4+\theta\left(1-(l_1/l_0)^2\right)-2\phi\left(1+l_1/l_0\right)\right)}.
\end{align*}
The inequalities \eqref{eq:f_cond_l0l1} are equivalent to
\begin{align*}
      &\theta \leqslant \frac{4}{\left(1+l_1/l_0\right)^2}, \\
      &\theta \leq \phi \leqslant \min \left[2\sqrt{\theta}-\frac{l_1}{l_0}\theta, \ \frac{4}{3\left(1+l_1/l_0\right)}+\frac{\theta}{3}\left(2-l_1/l_0\right)\right].
\end{align*}
if $\theta = \frac{4}{\left(1+l_1/l_0\right)\left(7+5l_1/l_0\right)}$ then maximum of $c(\phi)$ is attained at $\phi_* = \frac{4\left(2+l_1/l_0\right)}{\left(1+l_1/l_0\right)\left(7+5l_1/l_0\right)}$. If $\theta \neq \frac{4}{\left(1+l_1/l_0\right)\left(7+5l_1/l_0\right)}$ then the derivative $c^{,}_{\phi}$ has two roots:
\begin{align*}
&\phi_*^{+-}:= \\\\
&\frac{\frac{\theta^2 l_1}{2l_0}\left(1-\frac{l_1(1+l_1/l_0)}{2l_0}\right) + \theta\left(\theta+2+2l_1/l_0\right)\pm \left(1-\theta (1+l_1/l_0)^2/4\right)\sqrt{\theta\left[\theta+4+\theta\left(1-(l_1/l_0)^2\right)\right]}}{\frac{\theta}{4}\left(1+l_1/l_0\right)\left(5l_1/l_0+7\right)-1}
\end{align*}
If $\theta < \frac{4}{\left(1+l_1/l_0\right)\left(7+5l_1/l_0\right)}$ then maximum of $c(\phi)$ is attained at the bigger root:
\begin{align*}
    \theta \leqslant \phi_*^{-} \leqslant 2\sqrt{\theta}-\frac{l_1}{l_0}\theta.
\end{align*}
It's straightforward to check that the relevant restrictions are satisfied. The same way, for $\theta > \frac{4}{\left(1+l_1/l_0\right)\left(7+5l_1/l_0\right)}$, there are two roots. The smaller one is
\begin{align*}
    \theta \leqslant  \phi_*^{-} \leqslant \min \left[2\sqrt{\theta}-\frac{l_1}{l_0}\theta, \ \frac{4}{3\left(1+l_1/l_0\right)}+\frac{\theta}{3}\left(2-l_1/l_0\right)\right].
\end{align*}
The bigger root falls outside the admissible range of $\phi$:
\begin{align*}
    \min \left[2\sqrt{\theta}-\frac{l_1}{l_0}\theta, \ \frac{4}{3\left(1+l_1/l_0\right)}+\frac{\theta}{3}\left(2-l_1/l_0\right)\right] \leqslant \phi_*^{+}.
\end{align*}
In conclusion, the maximum of $c(\phi)$ is attained at $\phi_* = \phi_*^{-}$.

\endproof

Next we shall consider the case of $l_1 > l_0$. Unfortunately, the closed form expression for worst-case contraction seems to be unreasonably complicated. Therefore, we will just write the contraction expression for arbitrary $F$ and the range of $F$ that contains the maximal point. With that, it will be easier to estimate the worst-case contraction factor numerically, rather than using an overly complicated expression.

%Lemma for contraction factor upper bound: (L_0,L_1)-smooth convex functions (l0 < l1)
\begin{lem} \label{lem:smooth_contraction_gen_2} Suppose the conditions of Theorem~\ref{teo:general_FG_res} are satisfied. For set initial gradient norm $G$, bounded distance to the optimum $R$, the following upper bound holds:
\begin{align}
    c  &= \frac{f(x^{k+1})-f(x_*)}{F} \notag \\
    &\leqslant \sqrt{\frac{\alpha}{\beta}}\left(\frac{\sqrt{\alpha}-\sqrt{\beta}}{l_1/l_0-1}\right)^2 + \frac{\gamma}{l_1/l_0-1}\left(\frac{l_1}{l_0} - \sqrt{\frac{\alpha}{\beta}}\right), 
\label{eq:gen_contr_bound_2}
\end{align}
where 
\begin{align*}
    &\theta := \left(\frac{G}{l_1R}\right)^2 \leqslant \frac{4}{\left(1+l_1/l_0\right)^2}, \\
    &\alpha := \frac{1}{\phi_*}-\frac{\theta}{\phi_*}\left(\frac{l_1}{l_0}-1\right) - 1, \\
    &\beta := \frac{1}{\phi_*} - \frac{l_1}{l_0}, \\
    &\gamma = 1 - \frac{\theta}{\phi_*},
\end{align*}
and the point of maximum $\phi_*$ lies in the following range:
\begin{align*}
    \theta \leqslant \phi_* \leqslant \min \left[2\sqrt{\theta} - \frac{l_1}{l_0}\theta, \ \frac{l_1/l_0+1}{\left(l_1/l_0\right)^2+l_1/l_0+1}\left(1+\frac{\theta}{l_1/l_0+1}\right) \right].
\end{align*}
\end{lem}

\proof

Consider the upper bound for $f(x^{k+1}) -f(x_*)$ from the statement of Theorem~\ref{teo:general_FG_res_1}. First assume $F > \frac{l_1/l_0+1}{\left(l_1/l_0\right)^2+l_1/l_0+1}\left[\frac{l_1R^2}{2} + \frac{G^2}{2l_1}\frac{1}{l_1/l_0+1}\right]$. In this case which we call the second mode, the contraction expression is
\begin{align}
    \frac{f(x^{k+1}) - f(x_*)}{F} \leqslant \frac{1}{F}\left[\frac{l_0R^2}{2} - \frac{G^2}{2l_1}\left(1-l_0/l_1\right)\right] -\frac{l_0}{l_1}.
\end{align}
From the expression above, the bound for contraction is a decreasing function of $F$ which approaches maximum at $F = \frac{l_1/l_0+1}{\left(l_1/l_0\right)^2+l_1/l_0+1}\left[\frac{l_1R^2}{2} + \frac{G^2}{2l_1}\frac{1}{l_1/l_0+1}\right]$. Therefore we only need to consider $F \leqslant \frac{l_1/l_0+1}{\left(l_1/l_0\right)^2+l_1/l_0+1}\left[\frac{l_1R^2}{2} + \frac{G^2}{2l_1}\frac{1}{l_1/l_0+1}\right]$. By dividing the \eqref{eq:g_val_bound_2} by $F$ and denoting 
\begin{align*}
    &\theta \colon= \left(\frac{G}{l_1R}\right)^2, \ \phi \colon= \frac{2F}{l_1R^2},
\end{align*}
we get the required expression \eqref{eq:gen_contr_bound_2} for contraction factor. Per \eqref{eq:f_cond_l0l1_2}, the expression attains maximum inside the following boundaries:
\begin{align*}
    \theta \leqslant \phi_* \leqslant \min \left[2\sqrt{\theta} - \frac{l_1}{l_0}\theta, \ \frac{l_1/l_0+1}{\left(l_1/l_0\right)^2+l_1/l_0+1}\left(1+\frac{\theta}{l_1/l_0+1}\right) \right].
\end{align*}

\endproof

We use  results of Lemma~\ref{lem:smooth_contraction_gen_1} and Lemma~\ref{lem:smooth_contraction_gen_2} to improve one specific result for $(L_0,L_1)$-smooth convex functions. The Theorem 3.1. from \cite{lobanov2025linearconvergencerateconvex} states that if  $\|\nabla f(x^k)\|_2 \geqslant \frac{L_0}{L_1}$ then for $(L_0,L_1)$-GD with $\eta =1$:
\begin{align}
    f(x^{k+1}) - f(x_*) \leqslant \left(1-\frac{1}{4L_1R} \right) \left(f(x^{k}) - f(x_*) \right).
\label{eq:old_result}
\end{align}
Previously in Section~\ref{sec:construction} we set $l_0, \ l_1$ as
\begin{align*}
    &l_0 :=  L_0\exp{(L_1R)},\\
    &l_1 := \frac
    {\left(L_0+L_1G\right)}{\nu},
\end{align*}
where $\nu$ is a constant such that $\nu = \exp{(-\nu)}$, can be estimated as $\nu \approx 0,567143$. \par
\begin{figure}[!ht]
\centering
\includegraphics{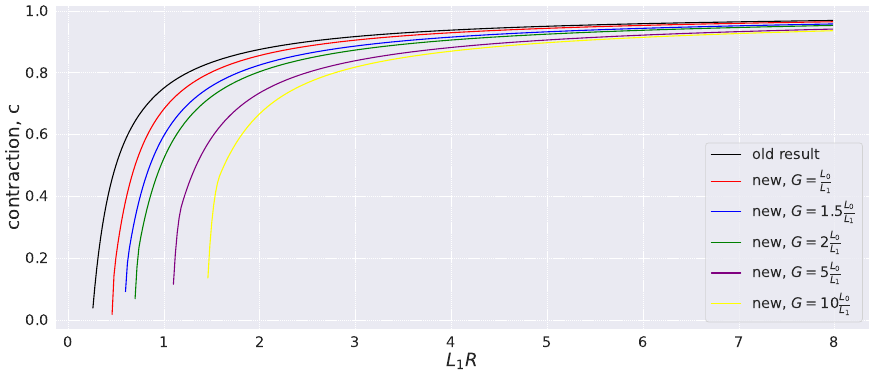} 
\caption{New contraction factor guarantees for several initial gradient's norms $\|\nabla f(x^0)\| = G$. The curves can not be extended to even smaller $L_1R$ because those values are outside the feasible range (smallest possible $L_1R$ is individual for each curve, Lemma~\ref{lem:feasibility_ness_gen1}). Graph $G = L_0/L_1$ is essentially, a novel estimate of contraction factor under gradient regime $\|\nabla f(x^0)\| > L_0/L_1$. } \label{fig:fig_3}
\end{figure}

\begin{figure}[!ht]
\centering
\includegraphics{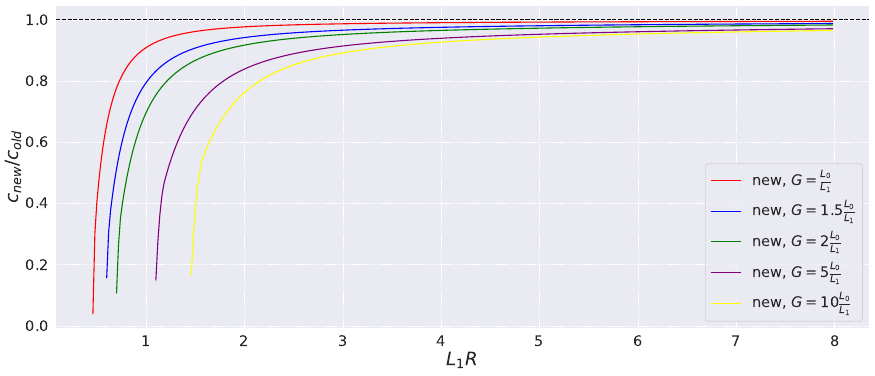} 
\caption{New contraction factor guarantees for several initial gradient's norms $\|\nabla f(x^0)\| = G$ divided by the old contraction factor.} \label{fig:fig_4}
\end{figure}

 
To apply our results we need to express $l_1/l_0$ and $\theta$ through $L_0, \ L_1, \ G$ and $R$:
\begin{align*}
    &\frac{l_1}{l_0} = \frac{1}{\nu} \left(1+\frac{L_1}{L_0}G\right)\exp{(-L_1R)}, \\
    &\theta = \left(\frac{G}{l_1R}\right)^2 = \left[\frac{\nu}{L_1R}\frac{\frac{L_1}{L_0}G}{1+\frac{L_1}{L_0}G}\right]^2.
\end{align*}
where $\theta$ should not exceed $\frac{4}{(1+l_1/l_0)^2}$. There are two possibilities: $l_1 \leqslant l_0$ which is equivalent to
\begin{align*}
    L_1R \geqslant \ln \left[\frac{1}{\nu}\left(1+\frac{L_1}{L_0}G\right)\right].
\end{align*}
we use contraction factor from Lemma~\ref{lem:smooth_contraction_gen_1} in this case. Otherwise, for $l_1 > l_0$ that is
\begin{align*}
    L_1R < \ln \left[\frac{1}{\nu}\left(1+\frac{L_1}{L_0}G\right)\right].
\end{align*}
we use the result of Lemma~\ref{lem:smooth_contraction_gen_2}. Unlike \eqref{eq:old_result} that does not differentiate as long as $G > \frac{L_0}{L_1}$, the novel result does. Therefore, we plot a family of functions that depend on $L_1R$ and are parameterized by $\frac{L_1}{L_0}G$. The results are depicted on Figs.~\ref{fig:fig_3},\ref{fig:fig_4}. Even though each curve corresponds to set $G = G_0$, the contraction guarantee is valid for all $G \geqslant G_0$ because of monotone dependency on $G$.
Thus, the curves $G = L_0/L_1$ on both figures correspond to the novel guarantee for the regime $G \geqslant L_0/L_1$. The improvement is considerable, albeit for smaller $L_1R$. As $L_1R$ grows, the difference becomes less and less significant. On the other hand, it is evident that the state-of-the-art result lacks nuance and guarantee can be improved significantly if the gradient's norm is taken into account. In conclusion, we shall note that this section's main purpose is to merely demonstrate the potential of our proposed approach and we did not aim for many immediate practical results here.

%The results are depicted on Figs.~\ref{fig:fig_1}--\ref{fig:fig_4}. The first two figures (Figs.~\ref{fig:fig_3},\ref{fig:fig_2}) directly compare the novel guarantee for contraction factor with the state-of-the-art guarantee. The improvement is considerable, albeit for smaller $L_1R$. As $L_1R$ grows, the difference becomes less and less significant. The last three figures (Figs.~\ref{fig:fig_3},\ref{fig:fig_4}) depict new contraction factor guarantees for different set $G$, as well as the old result that is valid for all $G > \frac{L_0}{L_1}$. From the plotted figures, it becomes evident that  In conclusion, we shall note that this section's main purpose is to merely demonstrate the potential of our proposed approach and we did not aim for many immediate practical results here.

%The upper bound from Lemma~\ref{lem:smooth_contraction} is depicted on Fig.~\ref{fig:contraction_smooth_1}.

%\begin{figure}[!ht]
%\centering
%\includegraphics{fig_1.pdf} 
%\caption{Comparing old contraction factor result for $G>L_0/L_1$ (Eq.~\ref{eq:old_result}) with the novel one.} \label{fig:fig_1}
%\end{figure}

% \vspace{\floatsep}
%\begin{figure}[!ht]
%\centering
%\includegraphics{fig_2.pdf} 
%\caption{New contraction factor for $G>L_0/L_1$ divided by the old.} %\label{fig:fig_2}
%\end{figure}

% \vspace{\floatsep}

 \vspace{\floatsep}

 %\vspace{\floatsep}
%\begin{figure}[!ht]
%\centering
%\includegraphics{fig_5.pdf} 
%\caption{Contraction factor for asymmetric $(L_0,L_1)$-smooth convex functions if $G > L_0/L_1$} \label{fig:fig_5}
%\end{figure}

\clearpage

\section{Conclusion}
This paper presents an meta-approach of deriving upper bounds for settings involving generalized smoothness with the use of parameterized template solutions. The proposed templates are already quite versatile and can be used for a range of different generalized smoothness assumptions. With that said, there is a lot of space for generalizations. For instance, the value of the proposed templates could be increased considerably if the step size was included in \ref{problem:G_PEP} as separate parameter. The further generalization of the proposed templates and design of the new ones are left to future research. Overall, we believe that our proposed approach could become a valuable instrument for acquiring guarantees in settings with non-standard (or SDP-incompatible) assumptions. After developing sufficiently general templates, one might be able to perform a streamlined comprehensive analysis of first-order systems that would normally require very cumbersome problem-specific analysis. Unfortunately, the approach has considerable negative sides. First of all, the templates inherently lead to conservatism because the initial problem has to be relaxed to fit into a template. Secondly, the templates do not absolve of problem-specific theoretical analysis completely. That is because the problem needs to first be reformulated in a way that would allow to apply a template. And lastly, as it is demonstrated in this paper, obtaining a template solution is not a trivial task. And the difficulty grows dramatically due to increase of the number of parameters in the template.

\section{Acknowledgements}
I am profoundly grateful to Eduard Gorbunov (MBZUAI, Abu Dhabi) for advice, discussions, and mentorship. I thank Adrien Taylor (adrien.taylor@inria.fr) for fruitful discussions.

\bibliographystyle{plain}

\bibliography{references}

\newpage
\appendix

\section{Construction of PEP template based on smooth convex minimization} \label{appendix:smooth_template}

In this section of the appendix, we present a PEP template that is based on minimization of smooth convex functions. Unlike classic problem setup for fixed step GD where only initial distance is bounded, we consider a situation when the gradient and function value in the starting point are explicit parameters of the problem. Although we do not apply this template to improve existing bounds, we believe that it still presents good theoretical value due to its potential for generalization. In fact, this was the first template we constructed, and we got the template for $(L_0,L_1)$-smoothness by introducing new parameters into it. \par

First we have to determine admissible initial function value and gradient's norm.  Let $\|x^k-x_*\|_2 \leqslant R$, $f(x^k)-f(x_*)=F$ and $\|\nabla f(x^k)\|_2=G$. From $\|\nabla f(x^k)-\nabla f(x_*) \|_2 \leqslant \|x^k-x_*\|_2$ we get $G \leqslant RL$. Suppose $G \leqslant RL$ is fixed. With
\begin{align*}
    &f(x_*) \geqslant f(x^k) + \langle \nabla f(x^k), x_{*}-x^k\rangle + \frac{\|\nabla f(x^k)-\nabla f(x_*) \|_2^2}{2L},\\
    &f(x^k) \geqslant f(x_*) + \langle \nabla f(x_*), x^{k}-x_*\rangle + \frac{\|\nabla f(x^k)-\nabla f(x_*) \|_2^2}{2L}.
\end{align*}
it's easy to prove that $F$ has to satisfy
\begin{align*}
    \frac{G^2}{2L}\leqslant F \leqslant GR - \frac{G^2}{2L}.
\end{align*}
The above conditions for $R,G,F$ are also sufficient for existence of a smooth convex function that satisfies the initial conditions. 

%Lemma: necessary conditions for feasibility
\begin{lem} \label{lem:feasibility_base} Suppose values $R>0$, $G>0$ and $F>0$ satisfy 
\begin{align}
    &G \leqslant RL, \\
    &\frac{G^2}{2L}\leqslant F \leqslant GR - \frac{G^2}{2L}.
\label{eq:base_sufficiency}
\end{align}
then there exists $f \in \mathcal F_L$ such that $\|x^0-x_*\|_2 \leqslant R$, $\|\nabla f(x^0)\|^2_2 = G$ and $f(x^0)-f(x_*)=F$.
\end{lem}
\proof
Suppose $x_* = 0$ without losing generality. We propose a $1$-dimensional function
\begin{align*}
h(x) = 
    \begin{cases}
        &0 \quad \textit{if} \quad  |x| < p_0 = R - \frac{F}{G} - \frac{G}{2L}, \\
        &\frac{L|x-p_0|^2}{2} \quad \textit{if} \quad p_0 \leqslant |x| < p_1 = R - \frac{F}{G} + \frac{G}{2L}, \\
        &\frac{G^2}{2L} + G|x| \quad \textit{if} \quad p_1 \leqslant |x| \leqslant R.
    \end{cases}
\end{align*}
It's easy to see that this function $h \in \mathcal F_L$ , satisfies all the conditions and exists for any $R,G,F$ that satisfy \eqref{eq:base_sufficiency}.
\endproof

Unlike functional constraints from Problem~\ref{problem:G_PEP}, a system of pairwise inequalities from Theorem~\ref{teo:smooth_convex_interpolation} constitute actual interpolation conditions. Therefore, there is no need to consider feasibility of the whole (PEP) problem that has three points $x^k, \ x^{k+1}, \ x_*$. It suffices to show that there exists a smooth convex function satisfying $\|x^k-x_*\|_2 \leqslant R$, $f(x^k)-f(x_*)=F$ and $\|\nabla f(x^k)\|_2=G$ and the feasibility of the PEP template follows as a corollary. This is because interpolation conditions (i.e. conditions that are necessary and sufficient for interpolability) have a property of \textit{point-wise extensibility} \cite{rubbens2024constraint}. The idea is the following: suppose there is a system of sets $\{(x^i,g^i,f^i)\}_{i}$ that satisfies (pairwise) interpolation constraints. Then for any new $x^N$, there will necessarily exist such $g_N$ and $f_N$ that the new system $\{(x^i,g^i,f^i)\}_{i} \ \cup \ \{(x^N,g^N,f^N)\}$ also satisfies the interpolation conditions \cite{rubbens2024constraint}. The takeaway from this is the following: if functional constraints in the template are (pairwise) interpolation conditions then it is enough to show that a there exists a function that interpolates $x^k$ and $x_*$ and satisfies initial conditions. Now we introduce upper bounds for the error of function value when one step of GD is performed.  

\begin{teo} \label{teo:smooth_FG_res} Let $f \in \mathcal F_L$  with $f(x^k)-f(x_*)=F$, $\|\nabla f(x^k)\|_2=G$ and $\|x^k-x_*\|_2 \leqslant R$ ($R>0$, $G>0$, $F>0$ and satisfy the restrictions that were discussed earlier) then for an iterate generated by GD with $\eta = 1$ it holds that
if 
\begin{align}
    &G \leqslant RL, \label{eq:smooth_mode_1_restr} \\
    &\frac{G^2}{2L} \leqslant F \leqslant  \min \left [\frac{1}{3}\left(LR^2+\frac{G^2}{2L}\right), GR - \frac{G^2}{2L}\right]. \notag
\end{align}
then
\begin{align}
    f(x^{k+1}) - f(x_*) \leqslant \left(F-\frac{G^2}{2L}\right)\left(1 - \frac{1}{4}(1-\frac{G^2}{2LF})/(\frac{LR^2}{2F}-1)\right). 
\label{eq:smooth_bound_1}
\end{align}
However if
\begin{align}
    &\frac{RL}{2} \leqslant G \leqslant RL, \\
    &\frac{1}{3}\left(LR^2+\frac{G^2}{2L}\right) \leqslant F \leqslant GR - \frac{G^2}{2L} \notag
    \label{eq:smooth_mode_2_restr}
\end{align}
then
\begin{align}
    f(x^{k+1}) - f(x_*) \leqslant \frac{LR^2}{2} - F.
\end{align}
\end{teo}
\proof
To prove the Theorem's statement, we consider properties of $L$-smooth convex functions and the additional constraints for the initial iterate $x^k$. 
\begin{align*}
    &\lambda_{s0}: f(x_*) \geqslant f(x^k) + \langle \nabla f(x^k), x_{*}-x^k\rangle + \frac{\|\nabla f(x^k)-\nabla f(x_*) \|_2^2}{2L}, \\
    &\lambda_{s1}: f(x_*) \geqslant f(x^{k+1}) + \langle \nabla f(x^{k+1}),x_{*}-x^{k+1}\rangle + \frac{\|\nabla f(x^{k+1})-\nabla f(x_*) \|_2^2}{2L}, \\
    &\lambda_{01}: f(x^{k}) \geqslant f(x^{k+1}) + \langle \nabla f(x^{k+1}), x^{k}-x^{k+1} \rangle + \frac{\|\nabla f(x^{k+1})-\nabla f(x^k) \|_2^2}{2L}, \\
    & \lambda_R:R^2 - \|x^{k}-x_{*} \|_2^2 \geqslant 0, \\
    &\lambda_G: G^2 - \|\nabla f(x^k)\|_2^2 = 0, \\
    &\lambda_F: F-\left(f(x^k)-f(x_*)\right) = 0.
\end{align*}
We will be adding up the inequalities with the corresponding multipliers
\begin{align*}
    &0 \leqslant \lambda_{s1} := \left(1-\frac{G^2}{2LF}\right)/\left(\frac{LR^2}{F}-2\right), \\
    &0 \leqslant \lambda_{s0} := \lambda_{s1}^2, \\
    &0 \leqslant \lambda_{01} := 1 - \lambda_{s1}, \\
    &\lambda_F = 1 - \lambda_{s1} - \lambda_{s0}, \\
    &0 \leqslant \lambda_{R} = \frac{L}{2}\lambda_{s0}, \\
    &\lambda_{G} = -\frac{\lambda_{01}}{2L}.
\end{align*}
and substituting $x^{k+1} = x^{k} - \nabla f(x^k)/L$ leads to 
\begin{align*}
    f(x^{k+1}) - f(x_*) \leqslant \left(F-\frac{G^2}{2L}\right)\left(1 - \frac{1}{4}(1-\frac{G^2}{2LF})/(\frac{LR^2}{2F}-1)\right) - p^\top (S\otimes I_n) p,
\end{align*}
where $p:=((x^k-x_{*})^\top, \left(\nabla f(x^k)-\nabla f(x_*))^\top,(\nabla f(x^{k+1})-\nabla f(x_*))^\top \right)^\top$ and 
\begin{align*}
    &S:= \frac{1}{2} \begin{bmatrix*}
	\lambda_{s0}L  & -\lambda_{s0} & -\lambda_{s1} \\
	-\lambda_{s0} & \lambda_{s0}/L & \lambda_{s1}/L \\
	-\lambda_{s1}  & \lambda_{s1}/L & 1/L
    \end{bmatrix*}
\\
\end{align*}
is a positive semidefinite matrix. We also need to determine the range of parameters that ensures  $\lambda_{s1}, \lambda_{s0}, \lambda_{01}, \lambda_{R}  \geqslant 0$. Obviously $\lambda_{R}$ and $\lambda_{s0}$ are always non-negative. For $\lambda_{s1}$:
\begin{align*}
    \frac{G^2}{2L} \leqslant F < \frac{LR^2}{2}
\end{align*}
which also implies $G < RL$ and from $\lambda_{01} \geqslant 0$:
\begin{align*}
     F \leqslant \frac{1}{3}\left(LR^2+\frac{G^2}{2L}\right).
\end{align*}
Now consider a different set of multipliers:
\begin{align*}
    &\lambda_{s0} = \lambda_{s1} = 1, \quad
    \lambda_{01} = \lambda_{G} = 0 \\
    &\lambda_F = -1, \quad
    \lambda_{R} = \frac{L}{2}
\end{align*}
By adding up the constraints, we get 
\begin{align*}
    f(x^{k+1}) - f(x_*) \leqslant \frac{LR^2}{2} - F -  p^\top (S\otimes I_n) p,
\end{align*}
where positive semidefinite matrix $S$ is
\begin{align*}
    &S:= \frac{1}{2} \begin{bmatrix*}
	L  & -1 & -1 \\
	-1 & 1/L & 1/L \\
	-1  & 1/L & 1/L
    \end{bmatrix*}.
\\
\end{align*}
This upper bound is valid for the case $F \geqslant \frac{1}{3}\left(LR^2+\frac{G^2}{2L}\right)$. It's straghtforward to check that the bounds intersect at exactly $F = \frac{1}{3}\left(LR^2+\frac{G^2}{2L}\right)$. From the bound above (and known properties of smooth convex functions) we have $F \leqslant \frac{LR^2}{2}$. Thus, from $F \geqslant \frac{1}{3}\left(LR^2+\frac{G^2}{2L}\right)$ we get $G\leqslant LR$ as expected. It's straightforward to check that $GR - \frac{G^2}{2L}$ only exceeds $\frac{1}{3}\left(LR^2+\frac{G0^2}{2L}\right)$ for $G > \frac{RL}{2}$ therefore the second upper bound is only valid for such gradient's norms.

\endproof

We can also obtain upper bound for the case when only initial gradient norm $\|\nabla f(x^k)\|_2$ is fixed.

\begin{teo} \label{teo:smooth_G_res} Let $f$ be $L$-smooth convex function with $\|\nabla f(x^k)\|_2=G$ and $\|x^k-x_*\|_2 \leqslant R$ ($R > 0$, $G > 0$) then for an iterate generated by GD with $\eta = 1$ it holds that 
\begin{align}
    &f(x^{k+1}) - f(x_*) \leqslant \left(GR-\frac{3G^2}{2L}\right) \hspace{10pt} \text{if} \hspace{8pt} G < \frac{3-\sqrt{5}}{2} RL, \\
    &f(x^{k+1}) - f(x_*) \leqslant \frac{3-\sqrt{5}}{4}\left(LR^2-\frac{G^2}{L}\right)  \hspace{10pt} \text{if} \hspace{8pt} \frac{3-\sqrt{5}}{2} RL \leqslant G \leqslant RL.
\label{eq:smooth_G_res}
\end{align}
\end{teo}
\proof
By taking derivative of the right side of \eqref{eq:smooth_bound_1}:
\begin{align*}
    \left(F-\frac{G^2}{2L}\right)\left(1 - \frac{1}{4}(1-\frac{G^2}{2LF})/(\frac{LR^2}{2F}-1)\right), 
\end{align*}
we find that the largest value is attained with
\begin{align*}
    &F_{worst} = \frac{\sqrt{5}-1}{10}LR^2+\frac{\sqrt{5}}{10}\frac{G^2}{L} \hspace{10pt} \text{if} \hspace{8pt} G \geqslant \frac{3-\sqrt{5}}{2} RL, \\
    &F_{worst} = GR - \frac{G^2}{2L} \hspace{10pt} \text{otherwise}.
\end{align*}
Those $F_{worst}$ yield the statement of the Theorem. Also the bounds can be proven directly: consider the inequalities
\begin{align*}
    &\lambda_{s0}: f(x_*) \geqslant f(x^k) + \langle \nabla f(x^k),x_{*}-x^k \rangle + \frac{\|\nabla f(x^k)-\nabla f(x_*) \|_2^2}{2L}, \\
    &\lambda_{s1}: f(x_*) \geqslant f(x^{k+1}) + \langle \nabla f(x^{k+1}),x_{*}-x^{k+1} \rangle + \frac{\|\nabla f(x^{k+1})-\nabla f(x_*) \|_2^2}{2L}, \\
    &\lambda_{01}: f(x^{k}) \geqslant f(x^{k+1}) + \langle \nabla f(x^{k+1}),x^{k}-x^{k+1} \rangle + \frac{\|\nabla f(x^{k+1})-\nabla f(x^k) \|_2^2}{2L}, \\
    & \lambda_R:R^2 - \|x^{k}-x_{*} \|_2^2 \geqslant 0, \\
    &\lambda_G: G^2 - \|\nabla f(x^k) \|_2^2 = 0. \\
\end{align*}
The first case is $G \geqslant \alpha RL$ where we denote $\alpha:= \frac{3-\sqrt{5}}{2}$. By choosing the following multipliers 
\begin{align*}
    &\lambda_{s1} = 1-\alpha,\\
    &\lambda_{s0} = \lambda_{01} = \alpha, \\
    &\lambda_{R} = \frac{L}{2}\alpha, \\
    &\lambda_{G} = -\frac{\alpha}{2L},
\end{align*}
substituting $x^{k+1} = x^{k} - \nabla f(x^k)/L$ and taking a sum of inequalities, we obtain
\begin{align*}
    f(x^{k+1}) - f(x_*) \leqslant \frac{\alpha}{2}\left(LR^2-\frac{G^2}{L}\right) - p^\top (S\otimes I_n) p,
\end{align*}
where positive semidefinite matrix $S$ is
\begin{align*}
    &S:= \frac{1}{2} \begin{bmatrix*}
	-\alpha L  & -\alpha & -(1-\alpha) \\
	-\alpha & \alpha/L & (1-\alpha)/L \\
	-(1-\alpha)  & (1-\alpha)/L & 1/L
    \end{bmatrix*}.
\\
\end{align*}
The second case is $G < \alpha RL$. The multipliers are 
\begin{align*}
    &\lambda_{s0} = \lambda_{01} = \theta, \quad \lambda_{s1} = 1 - \theta, \\
    &\lambda_{R} = \frac{G}{2R}, \quad \lambda_{G} = \frac{R}{2G} - \frac{3}{2L},
\end{align*}
where we denoted $\theta := 1 - \frac{G}{LR-G}$. The resulting inequality is
\begin{align*}
    f(x^{k+1}) - f(x_*) \leqslant \left(GR-\frac{3G^2}{2L}\right) - p^\top (S\otimes I_n) p,
\end{align*}
where positive semidefinite matrix $S$ is
\begin{align*}
    &S:= \frac{1}{2} \begin{bmatrix*}
	G/R  & -\theta & -(1-\theta) \\
	-\theta & R/G-3/L+2\theta/L & (1-\theta)/L \\
	-(1-\theta)  & (1-\theta)/L & 1/L
    \end{bmatrix*}.
\\
\end{align*}
Thus, the proof is complete.

\endproof

If neither initial function value nor initial norm of the gradient are fixed, we get an expected bound \cite{taylor2017smooth} $f(x^{k+1})-f(x_*) \leqslant \frac{LR^2}{6}$ by maximizing the bound from Theorem~\ref{teo:smooth_G_res} over $G$ ($G_{worst} = \frac{LR}{3}$). Because the results of Theorems~\ref{teo:smooth_FG_res} and \ref{teo:smooth_G_res} exactly match numerical PEP experiments (see Appendix~\ref{appendix:pep_experiments}), we are inclined to believe that those are exact worst-case guarantees. We do not provide a worst function that leads to an exactly matching lower bound. With that said, we assume that such a function is not one dimensional because there is no such $1$-dimensional Huber's function that would exactly match the upper bound (for both Theorems~\ref{teo:smooth_FG_res} and \ref{teo:smooth_G_res}). \par 

We can proceed to apply the template to calculate the worst contraction factor of function value for set $G$.

%Lemma for contraction factor upper bound: smooth convex functions
\begin{lem} \label{lem:smooth_contraction} Suppose the conditions of Theorem~\ref{teo:smooth_FG_res} are satisfied. For set initial gradient norm $G$ (and bounded distance to the optimum $R$), the following upper bound holds:
\begin{align}
    c = \frac{f(x^{k+1})-f(x_*)}{f(x^{k})-f(x_*)} &= \frac{f(x^{k+1})-f(x_*)}{F} \notag \\
    &\leqslant 1 - \frac{\theta}{\phi_*} - \frac{(\phi_*-\theta)^2}{4\phi_*(1-\phi_*)},
\label{eq:smooth_contr_bound}
\end{align}
where $\theta := \left(\frac{G}{LR}\right)^2 \leqslant 1$ and  $\phi_*$ is equal to 
\begin{align*}
    &\phi_* = \frac{-\theta(\theta+4) + \sqrt{\theta(\theta+4)(1-\theta)^2}}{1-6\theta} \quad \text{if} \quad  \theta \neq 1/6, \\
    &\phi_* := 1/2 \quad \text{if} \quad \theta = 1/6.
\end{align*}
\end{lem}
\proof
Consider the upper bound for $f(x^{k+1}) -f(x_*)$ from the statement of Theorem~\ref{teo:smooth_FG_res}. First assume $G>\frac{LR}{2}$ and $F \geqslant \frac{1}{3}\left(LR^2+\frac{G^2}{2L}\right)$. In this case which we call the second mode, the contraction expression is
\begin{align}
    \frac{f(x^{k+1})-f(x_*)}{F} \leqslant \frac{LR^2}{2F} -1.
\end{align}
From the expression above, the bound for contraction is a decreasing function of $F$ which attains maximum at $F = \frac{1}{3}\left(LR^2+\frac{G^2}{2L}\right)$. Therefore we only need to consider $F \leqslant \frac{1}{3}\left(LR^2+\frac{G^2}{2L}\right)$:
\begin{align}
    \frac{f(x^{k+1}) -f(x_*)}{F} \leqslant \left(1-\frac{G^2}{2LF}\right)\left(1 - \frac{1}{4}(1-\frac{G^2}{2LF})/(\frac{LR^2}{2F}-1)\right). \notag
\end{align}
We denote
\begin{align*}
    &\theta \colon= \left(\frac{G}{LR}\right)^2, \ \phi \colon= \frac{2F}{LR^2}.
\end{align*}
And rewrite right side expression as a function of $\phi$:
\begin{align*}
    c(\phi) := 1 - \frac{\theta}{\phi} - \frac{(\phi-\theta)^2}{4\phi(1-\phi)}.
\end{align*}
The inequalities \eqref{eq:smooth_mode_1_restr} are equivalent to
\begin{align*}
      &\theta \leqslant 1, \\
      &\theta \leq \phi \leqslant \min \left[2\sqrt{\theta}-\theta, \ \frac{2}{3}+\frac{\theta}{3}\right].
\end{align*}
Derivative of $c(\phi)$:
\begin{align*}
    c^{,}_{\phi}  = \frac{(6\theta-1)\phi^2 - 2\theta(\theta+4)\phi+\theta(\theta+4)}{4(\phi -1)^2\phi^2}.
\end{align*}
If $\theta = 1/6$ then maximum of $c(\phi)$ is attained at $\phi_* = 1/2$ which leads to $c(1/2) = 5/9$. If $\theta \neq 1/6$ then $c^{,}_{\phi}$ has two roots:
\begin{align*}
    \frac{\theta(\theta+4) \pm \sqrt{\theta(\theta+4)(1-\theta)^2}}{6\theta-1}.
\end{align*}
If $\theta < 1/6$ then maximum of $c(\phi)$ is attained at the bigger root:
\begin{align*}
    \theta \leqslant \frac{-\theta(\theta+4) + \sqrt{\theta(\theta+4)(1-\theta)^2}}{1-6\theta} \leqslant 2\sqrt{\theta}-\theta.
\end{align*}
It's straightforward to check that the relevant restrictions are satisfied. The same way, for $\theta > 1/6$, there are two roots. The smaller one is
\begin{align*}
    \theta \leqslant \frac{\theta(\theta+4) - \sqrt{\theta(\theta+4)(1-\theta)^2}}{6\theta-1} \leqslant \min \left[2\sqrt{\theta}-\theta, \ \frac{2}{3}+\frac{\theta}{3}\right].
\end{align*}
the larger root falls outside the admissible range of $\phi$:
\begin{align*}
    \min \left[2\sqrt{\theta}-\theta, \ \frac{2}{3}+\frac{\theta}{3}\right] \leqslant \frac{\theta(\theta+4) + \sqrt{\theta(\theta+4)(1-\theta)^2}}{6\theta-1}.
\end{align*}
The maximum of $c(\phi)$ is attained at 
\begin{align*}
    \phi_* = \frac{\theta(4+\theta) - (1-\theta)\sqrt{\theta(4+\theta)}}{6\theta-1}.
\end{align*}

\endproof

It's straightforward to check that Lemma~\ref{lem:smooth_contraction_gen_1} is correctly reduced to Lemma~\ref{lem:smooth_contraction} if one lets $l_0 = l_1$ (exactly as it should).

\section{PEP experiments} \label{appendix:pep_experiments}

Throughout this paper we used numerical PEP solutions to guess multipliers for the constraints in order to attain theoretical proofs of upper bounds. More specifically, we employed this computer-assisted approach for Theorems~\ref{teo:smooth_FG_res},~\ref{teo:smooth_G_res},~\ref{teo:general_FG_res} and~\ref{teo:general_FG_res_1}. Consider the problem of upper-bounding the error of smooth convex function after one iteration of GD. The best worst-case upper bound can be expressed \cite{taylor2017smooth} as a solution of the following problem

\begin{align*} \label{problem:L_PEP}\tag{L-PEP}
    \und{\sup}{\{x^i,\nabla f(x^i), f(x^i)\}_{i\in I}} &f(x^{k+1}) - f(x_*) \quad \text{such that}  \\
    &I = \{k, k+1,*\},\\
    & f(x^i) \geqslant f(x
    _j) + \langle \nabla f(x^j),x^i-x^j \rangle + \frac{\|\nabla f(x^i) - \nabla f(x^j) \|_2^2}{L}, \hspace{4pt} i,j \in I,\\
    &x^{k+1} = x^k - \frac{\nabla f(x^k)}{L}\\
    &\{x_*,\nabla f(x_*), f(x_*)\} = \{\textbf{0}, \textbf{0}, 0\},\\
    & \|x^k-x_*\|_2 \leqslant R, \\
    & \|\nabla f(x^k)\|_2 = G, \\
    & f(x^k)-f_* = F.
\end{align*}

Using the standard approach, (\ref{problem:L_PEP}) can be rewritten as an SDP and then dualized \cite{goujaud2023fundamental}. Solving a dual problem is convenient because it would immediately yield a certificate of the upper bound via dual variables. We introduce for $i, j \in I$:
\begin{align*}
    &A_{ij}  \colon = 
\begin{bmatrix} 
\overline{x}_i \\ \overline{x}_j \\ \overline{g}_i \\ \overline{g}_j
\end{bmatrix}^\mathrm{T}
\!\!\!A
\begin{bmatrix} 
\overline{x}_i \\ \overline{x}_j \\ \overline{g}_i \\ \overline{g}_j
\end{bmatrix}
,\hspace{8pt}
A_R:= \begin{bmatrix*}
1  & 0 & 0  \\
0 & 0 & 0  \\
0  & 0 & 0
\end{bmatrix*},\hspace{8pt}
A_G:= \begin{bmatrix*}
0  & 0 & 0  \\
0 & 1 & 0  \\
0  & 0 & 0
\end{bmatrix*} \\\\
&m_{ij} \colon = \overline{f}_i - \overline{f}_j,
\end{align*}
where 
\begin{align*}
    A(L):= \frac{1}{2} \begin{bmatrix*}
	0  & 0 & 0 & -1 \\
	1 & 0 & 0 & 1 \\
	0  & 0 & -1/L & 1/L \\
    -1 & 1 & 1/L & -1/L
    \end{bmatrix*}.
\end{align*}
and basis vectors:
\begin{align*}
    &\overline{x}_* = [0,0,0], \hspace{8pt} \overline{x}_k = [1,0,0], \hspace{8pt} \overline{x}_{k+1} = [1,-1/L,0], \\
    &\overline{g}_* = [0,0,0], \hspace{8pt} \overline{g}_k = [0,1,0], \hspace{8pt} \overline{g}_{k+1} = [0,0,1], \\
    &\overline{f}_* = [0,0], \hspace{8pt} \overline{f}_k = [1,0], \hspace{8pt} \overline{f}_{k+1} = [0,1].
\end{align*}
With that, we present the dual problem for (\ref{problem:L_PEP}):
\begin{align*} \label{problem:dual_PEP}\tag{dual-PEP}
    \und{\inf}{\{\lambda_{ij}\}_{i,j\in I}, \lambda_R,\lambda_G,\lambda_F} &\lambda_RR^2 + \lambda_GG^2 + \lambda_FF 
    \quad \text{subject to}  \\
    &I = \{k, k+1,*\},\\
    &\lambda_{ij} \geqslant 0, \hspace{4pt} \forall i,j \in I, \\
    &\lambda_R \geqslant 0, \\
    & \sum_{i,j\in I} \lambda_{ij}A_{ij} - \lambda_R A_R - \lambda_G A_G \preccurlyeq 0, \\
    &m_{k+1,*} + \sum_{i,j\in I} \lambda_{ij}m_{ij} - \lambda_F m_{k,*} \leqslant 0.
\end{align*}

By repeatedly solving the problem above, we guessed the theoretical expressions for $\lambda_{ij}, \lambda_R,\lambda_G$ and $\lambda_F$ and derived proof of Theorem~\ref{teo:smooth_FG_res}. For Theorem~\ref{teo:smooth_G_res} we solved (\ref{problem:dual_PEP}) with $\lambda_F = 0$. 

To prove the bounds from Theorems~\ref{teo:general_FG_res},~\ref{teo:general_FG_res_1} we make several tweaks in the problem; for calculation of $A_{*,k}, A_{*,k+1}$ we use $A = A(l_0)$, $A = A(l_1)$ for $A_{k,k+1}, A_{k+1,k}$  and $A = A(l_2)$ otherwise. As for the basis vectors: $\overline{x}_{k+1} = [1,-1/l,0]$. The necessary conditions for feasibility of \ref{problem:G_PEP} presented in Lemma~\ref{lem:feasibility_ness_gen1} were numerically verified in a similar manner: the performance metric was replaced with either $\|\nabla f(x^k)\|_2$ or $f(x^k)-f(x_*)$. To find the lowest feasible $F$, maximization was changed to minimization in the non-dual problem, so the constraints in the dual SDP were also changed appropriately.

This concludes the description of the PEP numerical experiments. The numerical calculations were performed with the standard semidefinite solver MOSEK \cite{mosek}.

\end{document}